\documentclass[preprint,12pt]{elsarticle}

\usepackage{amssymb}
\usepackage{amsmath} 
\usepackage{graphicx}
\usepackage{hyperref}
\usepackage{booktabs}
\usepackage{listings}
\usepackage{xcolor}
\usepackage{ragged2e}
\usepackage{float}
\usepackage{subcaption}
\usepackage{algorithm} 
\usepackage{algpseudocode} 
\usepackage{caption}
\usepackage{setspace}
\usepackage[utf8]{inputenc}
\usepackage[modulo]{lineno}
\journal{}

\begin{document}

\begin{frontmatter}



\title{Two-dimensional seepage analysis using a polygonal cell-based smoothed finite element method}

\author[inst1]{Yang Yang}
\affiliation[inst1]{organization={PowerChina Kunming Engineering Corporation Limited},
            city={Kunming},
            postcode={650051}, 
            state={Yunnan},
            country={China}}

\author[inst2]{Mingjiao Yan}
\affiliation[inst2]{organization={College of Water Conservancy and Hydropower Engineering, Hohai University},
            city={Nanjing},
            postcode={210098}, 
            state={Jiangsu},
            country={China}}

\affiliation[inst3]{organization={College of Harbour, Coastal and Offshore Engineering},
            city={Nanjing},
            postcode={210098}, 
            state={Jiangsu},
            country={China}}

\affiliation[inst4]{organization={Yunnan Institute of Water Resources and Hydropower Research},
            city={Kunming},
            postcode={650228}, 
            state={Yunnan},
            country={China}}
\author[inst1]{Zongliang Zhang}
\author[inst3]{Yinpeng Yin}
\author[inst1]{Qiang Liu}
\author[inst4]{Youliang Li}

\begin{abstract}
This study develops a polygonal cell-based smoothed finite element method (CSFEM) for two-dimensional seepage analyses in porous media, covering steady-state, transient, and free-surface problems. Wachspress interpolation on convex polygonal elements is combined with cell-based gradient smoothing, so that element matrices are assembled using boundary integrals only, avoiding in-element derivatives and improving robustness on distorted and locally refined meshes. To improve efficiency, a solution-driven adaptive refinement strategy is employed to concentrate resolution near steep hydraulic gradients and evolving wet–dry interfaces. Free-surface seepage is handled by a fixed-mesh iterative scheme that updates the wetted region and boundary conditions to track the phreatic surface. Benchmark tests validate the formulation against analytical solutions and high-fidelity FEM references. In steady seepage examples, the proposed polygonal CSFEM reproduces linear hydraulic-head fields to near machine precision and yields smaller head errors than conventional FEM at the same characteristic mesh size. In transient problems, accurate head evolution and stable time responses are obtained, while adaptive refinement efficiently resolves localized high-gradient zones. For free-surface cases, the method captures the phreatic-surface profile and seepage-face development reliably without remeshing. The quadtree refinement and adaptivity provide substantial efficiency gains in degrees of freedom and runtime for a prescribed accuracy level.
\end{abstract}




\end{frontmatter}


\section{Introduction}
\label{sec:Introduce}

Seepage analysis plays a central role in many geotechnical and hydraulic engineering applications, such as the safety assessment of embankment dams \citep{kheiri2020numerical,kalateh2024review}, tunnels \citep{li2019study,chen2024seepage}, foundations \citep{xu2022seepage,yuan2021numerical}, and underground structures \citep{wang2021reliability,zhang2019buoyancy}. Accurate prediction of hydraulic head and pore-pressure distributions is essential for evaluating uplift forces, piping potential, and overall stability. The finite element method (FEM) remains one of the most widely used numerical tools for seepage problems due to its flexibility in handling heterogeneous materials and complex boundary conditions \citep{zienkiewicz2005finite,liu2022eighty}.

Despite its broad applicability, conventional FEM formulations based on triangular or quadrilateral elements face several limitations. Triangular meshes are straightforward to generate but may exhibit lower accuracy, whereas quadrilateral meshes often require more manual control \citep{sunAutomaticQuad2015,docampo2020regularization}. Standard isoparametric elements are also sensitive to mesh distortion, since the Jacobian must remain positive and element angles must satisfy geometric constraints \citep{liu2003fem,zienkiewicz2000fem}. Overall, these drawbacks make the conventional FEM less robust and less efficient for problems with complex geometries and highly distorted meshes.

To alleviate the above mesh-sensitivity issues while retaining the simplicity of low-order discretizations, the smoothed finite element method (SFEM) introduces a gradient smoothing operation over appropriately constructed smoothing domains, in which gradients are averaged and—via the divergence theorem—evaluated through boundary integrals in the physical space. This treatment reduces the dependence on the inverse Jacobian inherent in isoparametric mappings, thereby relaxing geometric restrictions and improving robustness on distorted meshes \citep{liuSmoothedFiniteElementMechanics2007,liuTheoreticalAspectsSmoothed2007,zengLiuSFEManOverview2018}. Existing studies have further shown that SFEM can achieve higher accuracy and faster energy-norm convergence than the conventional FEM \citep{liuSmoothedFiniteElementMechanics2007,yan2025fast}. Consequently, SFEM is less sensitive to element distortion and can accommodate elements of arbitrary geometry, which is particularly appealing for seepage-flow analyses in complex geometries where locally steep hydraulic gradients and evolving wet–dry interfaces can exacerbate mesh-quality issues.

Several SFEM variants have been actively developed in the past few years and are commonly distinguished by how the smoothing domains are constructed, including the node-based SFEM (NSFEM) \citep{zhouVolumetricLockingFree2022,lyuImplicitStabilizedNSFEM2024}, the edge-based SFEM (ESFEM) \citep{jiangAssessmentESFEM2024,choiBoundaryEnhancementSFEM2025}, the face-based SFEM (FSFEM) for three-dimensional meshes \citep{huangFaceBasedSFEManASISimulation2023,jiangStrengthMetroBogie2025}, and the cell-based SFEM (CSFEM) \citep{cuiHighOrderCSFEM2021,zhaoNSidedPolygonalCSFEM2024}. Among these, CSFEM is particularly attractive in practice because it stays close to the standard FEM workflow by subdividing each element into smoothing subcells and evaluating element contributions via boundary (surface/line) integration, while requiring only limited modifications to conventional FEM assembly and offering improved robustness for complex analyses.

While CSFEM improves robustness through gradient smoothing, the overall flexibility and mesh quality of the discretization remain crucial for seepage simulations in complex geometries, which has motivated the increasing use of polygonal discretizations to alleviate meshing constraints and enhance robustness on locally refined or highly irregular grids. For example, polygonal scaled boundary formulations have been implemented in commercial FE platforms and applied to steady-state and transient seepage analyses on polygon/quadtree meshes \citep{ye2021psbfemabaqus,yang2022geofluidsseepage}. Beyond scaled-boundary approaches, conservative variational schemes have been developed directly on general polygonal meshes for Darcy flow, including lowest-order weak Galerkin methods with optimal convergence and continuous normal fluxes \citep{liu2018wgdarcy}, as well as mimetic finite difference discretizations designed for diffusion/Darcy operators on unstructured polygonal meshes and nonmatching local refinements \citep{kuznetsov2004mfd}. More recently, virtual element discretizations have provided another flexible route for Darcy flow on polygonal meshes, with rigorous analysis and extensions to fractured porous media \citep{wu2025fastvemdarcy}. Overall, these studies indicate that polygonal elements can offer tangible advantages for seepage modelling; coupling such geometric flexibility with the CSFEM smoothing framework therefore provides a promising route to simultaneously enhance discretization versatility and numerical robustness in the present polygonal CSFEM approach.

This need for geometric flexibility and mesh-robust discretizations becomes even more pronounced in free-surface seepage, where an evolving wet--dry interface and localized steep hydraulic gradients lead to a nonlinear free-boundary problem with inequality constraints on partially saturated or outflow boundaries \citep{neuman1970femfree,bath1979nomeshiter,lacy1987freesurface}. Classical FEM treatments include iterative free-surface updating schemes and fixed-domain formulations such as the residual-flow procedure to avoid repeated remeshing \citep{neuman1970femfree,bath1979nomeshiter,lacy1987freesurface,desai1983residualflow}, while more recent strategies employ interface-capturing and moving-boundary techniques, including level-set and moving-mesh methods \citep{herreros2006levelset,frolkovic2012levelset,darbandi2007movingmesh}. Because the free surface and seepage face are solution-dependent and may exhibit gradient singularities near wet--dry transitions, locally refined meshes are often required around active-set switching regions and geometric singularities; accordingly, adaptive refinement remains an effective route to improve efficiency and robustness \citep{rank1986adaptive,boeriu2005adaptiveparallel,ashby2021adaptive}. These considerations motivate an adaptive extension of the present polygonal CSFEM framework for free-surface seepage, where refinement can be driven by solution features such as phreatic-surface curvature, seepage-face flux, and residual indicators while preserving the flexibility of quadtree--polygonal discretizations.

The objective of this work is to develop a polygonal CSFEM for seepage analysis in saturated porous media, using Wachspress interpolation on convex polygonal elements and gradient smoothing on triangular subcells so that smoothed hydraulic head gradients are obtained from boundary integrals only. Hybrid quadtree meshes are employed to handle hanging nodes without additional constraint equations, and a fixed-mesh iterative procedure is adopted for free-surface seepage. The proposed approach is validated through benchmark problems for steady-state, transient, and free-surface seepage. The remainder of this paper is organized as follows: Section~\ref{sec:CSFEM_seepage} presents the polygonal CSFEM formulation, Sections~\ref{sec:wachspress} and \ref{sec:free_surface} describe the Wachspress interpolation and free-surface strategy, Section~\ref{numerical_examples} reports numerical validations, and Section~\ref{Conclusions} concludes the work.

\section{Theoretical formulation of the polygonal CSFEM for seepage problems}
\label{sec:CSFEM_seepage}

We consider a porous medium occupying the domain $\Omega \subset \mathbb{R}^2$, bounded by $\Gamma=\Gamma_h\cup\Gamma_q$ with $\Gamma_h\cap\Gamma_q=\varnothing$, where $\Gamma_h$ and $\Gamma_q$ are the portions of the boundary on which the hydraulic head and the outward normal flux are prescribed, respectively. The transient seepage problem is described by Darcy's law and mass conservation, leading to the diffusion-type governing equation \cite{johari2018reliability}
\begin{equation}
\nabla \cdot (\mathbf{k}\nabla h) + p - S_s \frac{\partial h}{\partial t} = 0 
\quad \text{in } \Omega ,
\label{eq:gov-seep-CSFEM}
\end{equation}
where $\mathbf{k}$ is the hydraulic conductivity tensor, $h$ is the hydraulic head, $S_s$ denotes the specific storage coefficient, and $p$ represents the volumetric source or sink term.

The boundary conditions consist of prescribed-head and prescribed-flux boundaries
\begin{equation}
h = \bar{h}
\quad \text{on } \Gamma_h ,
\label{eq:bc-head-CSFEM}
\end{equation}
\begin{equation}
\mathbf{k}\nabla h \cdot \mathbf{n} = \bar{q}
\quad \text{on } \Gamma_q ,
\label{eq:bc-flux-CSFEM}
\end{equation}
where $\bar{h}$ is the prescribed hydraulic head on $\Gamma_h$, $\bar{q}$ is the prescribed outward normal flux on $\Gamma_q$, and $\mathbf{n}$ is the outward unit normal. For transient analyses, the initial condition is specified as
\begin{equation}
h(\mathbf{x},0)=h_0(\mathbf{x}) \quad \text{in } \Omega .
\label{eq:ic-head-CSFEM}
\end{equation}

\subsection{Galerkin weak form and spatial discretization}

Let $w$ be an admissible weighting function such that $w=0$ on $\Gamma_h$. Multiplying Eq.~\eqref{eq:gov-seep-CSFEM} by $w$ and integrating over $\Omega$ gives
\begin{equation}
\int_{\Omega} w \, \nabla \cdot (\mathbf{k}\nabla h)\, d\Omega
+ \int_{\Omega} w\, p \, d\Omega
- \int_{\Omega} w\, S_s\, \frac{\partial h}{\partial t}\, d\Omega = 0 .
\label{eq:weak-raw-CSFEM}
\end{equation}

Applying the divergence theorem to the first term yields
\begin{equation}
\int_{\Omega} w \, \nabla \cdot (\mathbf{k}\nabla h)\, d\Omega
= \int_{\Gamma} w\, \mathbf{k}\nabla h \cdot \mathbf{n} \, d\Gamma
- \int_{\Omega} (\nabla w)^T \mathbf{k}\nabla h \, d\Omega .
\label{eq:div-CSFEM}
\end{equation}

With $w=0$ on $\Gamma_h$ and using the Neumann condition on $\Gamma_q$, the Galerkin weak form can be written as
\begin{equation}
\int_{\Omega} (\nabla w)^T \mathbf{k}\nabla h \, d\Omega
+ \int_{\Omega} w\, S_s\, \frac{\partial h}{\partial t}\, d\Omega
= \int_{\Omega} w\, p \, d\Omega
+ \int_{\Gamma_q} w\, \bar{q} \, d\Gamma .
\label{eq:weak-form-CSFEM}
\end{equation}

Let $\Omega$ be discretized into polygonal elements $\{\Omega_e\}$. The trial and test functions are approximated by
\begin{equation}
h^h(\mathbf{x},t)=\sum_{I=1}^{n_e} N_I(\mathbf{x})\, h_I(t), 
\qquad 
w^h(\mathbf{x})=\sum_{J=1}^{n_e} N_J(\mathbf{x})\, \delta h_J ,
\label{eq:trial-test-CSFEM}
\end{equation}
where $N_I$ are polygonal shape functions (e.g., Wachspress interpolants on convex polygons), $n_e$ is the number of nodes of $\Omega_e$, and $\delta h_J$ are arbitrary nodal variations satisfying essential boundary conditions. Substituting Eq.~\eqref{eq:trial-test-CSFEM} into Eq.~\eqref{eq:weak-form-CSFEM} and assembling over all elements leads to the standard semi-discrete system
\begin{equation}
\mathbf{K}\mathbf{H}
+ \mathbf{M}\dot{\mathbf{H}}
= \mathbf{F} ,
\label{eq:semi-CSFEM}
\end{equation}
where $\mathbf{H}$ collects the global nodal heads, $\dot{\mathbf{H}}$ is its time derivative, and $\mathbf{F}$ includes source and Neumann contributions.

\subsection{Smoothed gradient over polygonal smoothing cells}

For each polygonal element $\Omega_e$, a cell-based smoothing strategy is employed. The element is subdivided into $n_c$ non-overlapping smoothing cells $\{\Omega_C\}_{C=1}^{n_c}$ such that $\Omega_e=\cup_{C=1}^{n_c}\Omega_C$, with boundaries $\Gamma_C$ and areas $A_C$. A piecewise-constant smoothing function is defined by \cite{liu2007smoothed}
\begin{equation}
\Phi(\mathbf{x}-\mathbf{x}_C)=
\begin{cases}
1/A_C, & \mathbf{x}\in \Omega_C,\\
0, & \mathbf{x}\notin \Omega_C,
\end{cases}
\label{eq:smoothing-func-CSFEM}
\end{equation}
where $\mathbf{x}_C$ is a representative point of the smoothing cell. The smoothed gradient of the hydraulic head is defined as the weighted-average of the compatible gradient
\begin{equation}
\widetilde{\nabla} h (\mathbf{x}_C)
= \int_{\Omega_C} \nabla h(\mathbf{x}) \, \Phi(\mathbf{x}-\mathbf{x}_C)\, d\Omega .
\label{eq:smooth-grad-CSFEM}
\end{equation}

Applying integration by parts to Eq.~\eqref{eq:smooth-grad-CSFEM} gives
\begin{equation}
\widetilde{\nabla} h (\mathbf{x}_C)
= \int_{\Gamma_C} h(\mathbf{x})\, \mathbf{n}(\mathbf{x})\,\Phi(\mathbf{x}-\mathbf{x}_C)\, d\Gamma
- \int_{\Omega_C} h(\mathbf{x})\, \nabla\Phi(\mathbf{x}-\mathbf{x}_C)\, d\Omega .
\label{eq:smooth-parts-CSFEM}
\end{equation}

Substituting Eq.~\eqref{eq:smoothing-func-CSFEM} into Eq.~\eqref{eq:smooth-parts-CSFEM}, and noting that $\nabla\Phi(\mathbf{x}-\mathbf{x}_C)=\mathbf{0}$ inside $\Omega_C$ for the constant smoothing function, we obtain the boundary-integral form
\begin{equation}
\widetilde{\nabla} h (\mathbf{x}_C)
= \frac{1}{A_C} \int_{\Gamma_C} h(\mathbf{x})\, \mathbf{n}(\mathbf{x}) \, d\Gamma ,
\label{eq:smooth-boundary-CSFEM}
\end{equation}
where $\mathbf{x}_C$ is a representative point of the smoothing cell, $\Gamma_C=\partial\Omega_C$ is its boundary, and $\mathbf{n}(\mathbf{x})$ is the outward unit normal on $\Gamma_C$. This boundary-integral representation involves only values of $h$ on $\Gamma_C$, so no in-cell evaluation of shape-function derivatives $\nabla N_I$ is required.

Using the interpolation $h^h(\mathbf{x})=\sum_{I=1}^{n_e} N_I(\mathbf{x}) h_I$, Eq.~\eqref{eq:smooth-boundary-CSFEM} yields
\begin{equation}
\widetilde{\nabla} h^h (\mathbf{x}_C)
= \sum_{I=1}^{n_e} \widetilde{\mathbf{B}}_I^{(C)}\, h_I ,
\label{eq:smooth-grad-shape-CSFEM}
\end{equation}
where the smoothed gradient operator is
\begin{equation}
\widetilde{\mathbf{B}}_I^{(C)}
= \frac{1}{A_C} \int_{\Gamma_C} N_I(\mathbf{x})\, \mathbf{n}(\mathbf{x}) \, d\Gamma .
\label{eq:BI-CSFEM}
\end{equation}

For straight boundary segments $e=1,\dots,m_C$ of $\Gamma_C$ with length $L_e$ and constant outward normal $\mathbf{n}_e$, one obtains
\begin{equation}
\widetilde{\mathbf{B}}_I^{(C)}
= \frac{1}{A_C} \sum_{e=1}^{m_C} 
\left( \int_{e} N_I(\mathbf{x})\, d\Gamma \right)\mathbf{n}_e .
\label{eq:BI-edge-sum-CSFEM}
\end{equation}

The edge integral can be evaluated by one-dimensional Gauss quadrature
\begin{equation}
\int_{e} N_I(\mathbf{x}) \, d\Gamma
\approx 
\sum_{g=1}^{n_g} N_I(\mathbf{x}_g)\, w_g \, L_e ,
\label{eq:edge-gauss-CSFEM}
\end{equation}
where $\{\mathbf{x}_g,w_g\}$ are Gauss points and weights on the edge.

\subsection{Element matrices with smoothed gradients}

In the Galerkin setting, both the trial and test function gradients are replaced by their smoothed counterparts over each smoothing cell. Defining the cell-wise smoothed gradient matrix
\begin{equation}
\widetilde{\mathbf{B}}^{(C)}
= 
\begin{bmatrix}
\widetilde{\mathbf{B}}_1^{(C)} &
\widetilde{\mathbf{B}}_2^{(C)} &
\cdots &
\widetilde{\mathbf{B}}_{n_e}^{(C)}
\end{bmatrix},
\label{eq:B-CSFEM}
\end{equation}
the smoothed element stiffness matrix is obtained as
\begin{equation}
\mathbf{K}_e
= \sum_{C=1}^{n_c} 
A_C \, 
\widetilde{\mathbf{B}}^{(C)T} 
\mathbf{k} \,
\widetilde{\mathbf{B}}^{(C)} .
\label{eq:Ke-CSFEM}
\end{equation}

The capacity matrix is computed using the standard interpolation of $h$ as
\begin{equation}
\mathbf{M}_e
= \int_{\Omega_e} N^T S_s N \, d\Omega
= \sum_{C=1}^{n_c}
\int_{\Omega_C} N^T S_s N \, d\Omega ,
\label{eq:Me-CSFEM}
\end{equation}
where the subcell integral may be evaluated by standard area quadrature on $\Omega_C$.

The elemental source and Neumann vectors are
\begin{equation}
\mathbf{f}_{p,e} 
= \int_{\Omega_e} N^T p \, d\Omega ,
\label{eq:fp-CSFEM}
\end{equation}
\begin{equation}
\mathbf{f}_{q,e} 
= \int_{\Gamma_{q,e}} N^T \bar{q} \, d\Gamma .
\label{eq:fq-CSFEM}
\end{equation}

After assembly over all elements and enforcement of essential boundary conditions on $\Gamma_h$, the global system takes the form of Eq.~\eqref{eq:semi-CSFEM} with $\mathbf{F}=\mathbf{f}_p+\mathbf{f}_q$.

\subsection{Time integration}

Let the time step be $\Delta t$ and denote the solution at $t^n$ by $\mathbf{H}^n$. Using the generalized $\theta$-method \cite{stuart1991dynamics}, the time derivative is approximated by
\begin{equation}
\dot{\mathbf{H}}^{\,n+\theta} \approx \frac{\mathbf{H}^{n+1}-\mathbf{H}^{n}}{\Delta t},
\qquad 
\mathbf{H}^{n+\theta} = (1-\theta)\mathbf{H}^n+\theta \mathbf{H}^{n+1},
\label{eq:theta-method-CSFEM}
\end{equation}
with $\theta\in(0,1]$. Evaluating $\mathbf{K}$, $\mathbf{M}$, and $\mathbf{F}$ at $t^{n+\theta}$ and applying the $\theta$-method to Eq.~\eqref{eq:semi-CSFEM} yields
\begin{equation}
\left(\frac{1}{\Delta t}\mathbf{M}^{n+\theta}+\theta\,\mathbf{K}^{n+\theta}\right)\mathbf{H}^{n+1}
=
\mathbf{F}^{n+\theta}
+
\left(\frac{1}{\Delta t}\mathbf{M}^{n+\theta}-(1-\theta)\,\mathbf{K}^{n+\theta}\right)\mathbf{H}^{n}.
\label{eq:theta-scheme-CSFEM}
\end{equation}

In particular, for the fully implicit backward Euler scheme ($\theta=1$), Eq.~\eqref{eq:theta-scheme-CSFEM} reduces to
\begin{equation}
\left( \mathbf{K}^{n+1} + \frac{1}{\Delta t}\mathbf{M}^{n+1} \right)\mathbf{H}^{n+1}
=
\mathbf{F}^{n+1} + \frac{1}{\Delta t}\mathbf{M}^{n+1}\mathbf{H}^{n}.
\label{eq:backward-euler-CSFEM}
\end{equation}

\section{Polygonal interpolation using Wachspress basis functions}
\label{sec:wachspress}

\begin{figure}[H]
  \centering
  \includegraphics[width=0.65\textwidth]{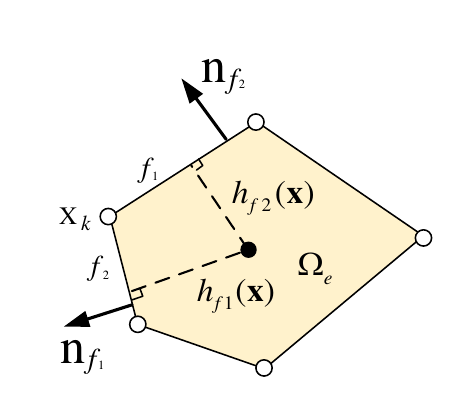}
  \caption{Illustration of Wachspress shape functions on a convex polygon.}
  \label{fig:shape_function}
\end{figure}

Consider a convex $n$-gon $\Omega_e$ with ordered vertices $\{\mathbf{x}_k\}_{k=1}^{n}$ in the physical plane. Fig.~\ref{fig:shape_function} provides a schematic illustration of the Wachspress shape functions on $\Omega_e$ \citep{warren2007barycentric}
. Let $f_1$ and $f_2$ denote the two edges adjacent to vertex $k$ and let $\mathbf{n}_{f_i}$ be the corresponding outward unit normals. For a point $\mathbf{x}\in\Omega_e$, define the signed distance to the supporting line of edge $f_i$ as
\begin{equation}
h_{f_i}(\mathbf{x})
= \big(\mathbf{x}-\mathbf{x}_{f_i}\big)\cdot \mathbf{n}_{f_i},
\label{eq:wach-hfi}
\end{equation}
where $\mathbf{x}_{f_i}$ is any point on edge $f_i$ (e.g., one of its endpoints).
With this convention, $h_{f_i}(\mathbf{x})>0$ for $\mathbf{x}$ inside $\Omega_e$.

The Wachspress weight associated with vertex $k$ is then defined as
\begin{equation}
w_k(\mathbf{x})
= \frac{\det\!\big(\mathbf{n}_{f_1},\,\mathbf{n}_{f_2}\big)}
       {h_{f_1}(\mathbf{x})\, h_{f_2}(\mathbf{x})},
\label{eq:wach-wk}
\end{equation}
where $\det(\mathbf{a},\mathbf{b})$ denotes the $2\times2$ determinant. For a fixed polygon, $\det(\mathbf{n}_{f_1},\mathbf{n}_{f_2})$ is a constant that depends only on the local edge geometry at vertex $k$.

The Wachspress shape function associated with vertex $k$ is given by
\begin{equation}
N_k(\mathbf{x})
= \frac{w_k(\mathbf{x})}
       {\displaystyle\sum_{j=1}^{n} w_j(\mathbf{x})}.
\label{eq:wach-Nk}
\end{equation}

Using these basis functions, a scalar field $h(\mathbf{x})$ is interpolated as
\begin{equation}
h(\mathbf{x})
= \sum_{k=1}^{n} N_k(\mathbf{x})\, h_k ,
\label{eq:wach-interp}
\end{equation}
where $h_k$ is the nodal value at vertex $k$. 

For convex polygons, the Wachspress basis is nonnegative, $C^0$-continuous, and linearly complete, and it satisfies the partition-of-unity property
\begin{equation}
\sum_{k=1}^{n} N_k(\mathbf{x}) = 1
\quad \text{for all } \mathbf{x} \in \Omega_e .
\label{eq:wach-pu}
\end{equation}

In the polygonal CSFEM, Wachspress functions enter the formulation only through boundary integrals of the form
\begin{equation}
\widetilde{\mathbf{B}}_I^{(C)}
= \frac{1}{A_C} \int_{\Gamma_C} N_I(\mathbf{x})\,\mathbf{n}(\mathbf{x})\, d\Gamma ,
\label{eq:wach-CSFEM-BI}
\end{equation}
and their edge-wise representation
\begin{equation}
\widetilde{\mathbf{B}}_I^{(C)}
= \frac{1}{A_C} \sum_{e=1}^{m_C}
\left(\int_{e} N_I(\mathbf{x})\, d\Gamma\right)\mathbf{n}_e ,
\label{eq:wach-CSFEM-BI-edge}
\end{equation}
which involve only $N_I$ but not $\nabla N_I$. Since $N_I$ is evaluated along straight edges, the integrand is smooth on each edge,
and accurate results can be obtained using a small number of Gauss points.

\section{Free-surface solutions}
\label{sec:free_surface}

For steady seepage with a phreatic surface, the computational domain $\Omega$ is separated by an unknown free surface into a wetted region $\Omega_w$ and a dry region $\Omega_d$. For simplicity, an isotropic conductivity $k$ is assumed in the free-surface examples below. The steady governing equation in the wetted region is written as
\begin{equation}
\nabla \cdot \left( k \nabla h \right) = 0 \quad \text{in } \Omega_w.
\label{eq:fs_gov}
\end{equation}

\begin{figure}[H]
  \centering
  \includegraphics[width=0.8\textwidth]{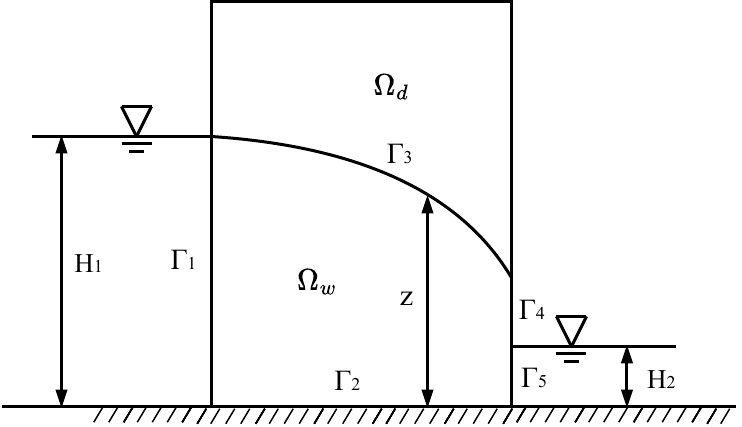}
  \caption{Geometry and boundary conditions for the free-surface seepage model.}
  \label{fig:Free_surfaces}
\end{figure}

Let the boundary of the seepage domain be decomposed as
$\Gamma = \Gamma_1 \cup \Gamma_2 \cup \Gamma_3 \cup \Gamma_4 \cup \Gamma_5$,
where $\Gamma_1$ is the upstream prescribed-head boundary, $\Gamma_2$ is the impervious boundary, $\Gamma_3$ is the unknown free surface, $\Gamma_4$ is the downstream seepage face, and $\Gamma_5$ is the downstream prescribed-head boundary representing the tailwater level. The boundary conditions are given by
\begin{equation}
h = H_1 \quad \text{on } \Gamma_1 ,
\label{eq:fs_bc_S1}
\end{equation}
\begin{equation}
h = H_2 \quad \text{on } \Gamma_5 ,
\label{eq:fs_bc_S5}
\end{equation}
\begin{equation}
-k \nabla h \cdot \mathbf{n} = 0 \quad \text{on } \Gamma_2 ,
\label{eq:fs_bc_S2}
\end{equation}
and the zero pressure-head condition on both the free surface and the seepage face,
\begin{equation}
h = z \quad \text{on } \Gamma_3 \cup \Gamma_4 ,
\label{eq:fs_bc_S34}
\end{equation}
where $z$ denotes the elevation coordinate. The free surface $\Gamma_3$ is not known a priori and is determined as part of the solution.

In this work, a fixed-mesh approach is adopted to determine the free surface $\Gamma_3$ by iteratively updating the wetted region while keeping the mesh topology unchanged. This strategy is well suited for the polygonal CSFEM because the formulation relies on boundary integrals and can robustly accommodate quadtree-based local refinements and hanging nodes through hybrid polygonal elements. The complete fixed-mesh iteration procedure, with an optional adaptive refinement loop, is summarized in Algorithm~\ref{alg:fs_adapt}.

At iteration $t$, a permeability field is assigned on the fixed mesh as
\begin{equation}
k^{(t)}(\mathbf{x})=
\begin{cases}
k, & \mathbf{x}\in \Omega_w^{(t)},\\
\alpha k, & \mathbf{x}\in \Omega_d^{(t)},
\end{cases}
\quad \text{with } \alpha \ll 1 ,
\label{eq:fs_k_update}
\end{equation}
where $\alpha$ is a small parameter and is taken as $\alpha=10^{-3}$ throughout this work. After solving for $h^{(t)}$, the pressure head is computed as $\psi^{(t)}(\mathbf{x})=h^{(t)}(\mathbf{x})-z$. 

The saturated region is then updated by
\begin{equation}
\Omega_w^{(t+1)} = \left\{\mathbf{x}\in\Omega \, \big| \, \psi^{(t)}(\mathbf{x}) \ge 0 \right\}, \qquad \Omega_d^{(t+1)}=\Omega\setminus \Omega_w^{(t+1)} .
\label{eq:fs_sat_update}
\end{equation}

The seepage face on the downstream vertical boundary $\Gamma_4$ is updated as the portion where $\psi^{(t)} \ge 0$ (below the intersection point), on which the atmospheric condition Eq.~\eqref{eq:fs_bc_S34} is enforced, while the remaining part of $\Gamma_4$ above the intersection point is treated as impermeable using Eq.~\eqref{eq:fs_bc_S2}. The prescribed-head condition on $\Gamma_5$ remains unchanged during the iteration.

The iteration is terminated when the change of the free-surface location becomes sufficiently small. In implementation, this can be monitored by the displacement of the downstream intersection point $x_o$:
\begin{equation}
\left|x_o^{(t+1)}-x_o^{(t)}\right| < \varepsilon_x .
\label{eq:fs_stop}
\end{equation}

\subsection{Adaptive refinement indicator}
The adaptive refinement step is optional and aims to improve resolution around the active-set switching region and steep hydraulic gradients. Let $K$ denote an element in the current mesh $\mathcal{T}_h$. We first define a free-surface band indicator
\begin{equation}
I_K =
\begin{cases}
1, & \min\limits_{\mathbf{x}\in \mathcal{S}(K)}\psi(\mathbf{x}) \le 0 \le \max\limits_{\mathbf{x}\in \mathcal{S}(K)}\psi(\mathbf{x}),\\
0, & \text{otherwise},
\end{cases}
\label{eq:band_indicator}
\end{equation}
where $\mathcal{S}(K)$ denotes a set of sampling points on $K$. This indicator detects elements intersected by the zero level set $\psi=0$.

To prioritize refinement within this band, a gradient indicator is further defined as
\begin{equation}
G_K = \left\|\nabla h\right\|_{0,K}
= \left(\int_{K}\left|\nabla h\right|^2\,\mathrm{d}\Omega\right)^{1/2},
\label{eq:grad_indicator}
\end{equation}
where $\nabla h$ is evaluated from the smoothed hydraulic-gradient operator on $K$. The refinement indicator is then taken as
\begin{equation}
\eta_K = I_K\,G_K,
\label{eq:eta_indicator}
\end{equation}
so that refinement is restricted to a narrow band around the free surface while emphasizing elements with large hydraulic gradients.

\subsection{Marking strategy and stopping criterion}
Given $\eta_K$ on all $K\in\mathcal{T}_h$, a bulk-chasing (D\"orfler) marking strategy is adopted \citep{Dorfler1996,Morin2002}: we select a minimal set of elements $\mathcal{M}\subset \mathcal{T}_h$ such that
\begin{equation}
\sum_{K\in\mathcal{M}}\eta_K^2 \ge \theta \sum_{K\in\mathcal{T}_h}\eta_K^2,
\label{eq:dorfler_marking}
\end{equation}
with $\theta\in(0,1)$, and refine all elements in $\mathcal{M}$ while enforcing a 2:1 balance in the quadtree mesh. The optional adaptivity loop is terminated when the fixed-mesh iteration satisfies Eq.~\eqref{eq:fs_stop} and the global indicator norm
\begin{equation}
\|\eta\|_{\ell^2(\mathcal{T}_h)} =
\left(\sum_{K\in\mathcal{T}_h}\eta_K^2\right)^{1/2}
\label{eq:eta_norm}
\end{equation}
falls below a prescribed tolerance. In practice, the adaptivity loop may also be stopped when the marked set $\mathcal{M}$ becomes empty.

\begin{algorithm}[H]
\caption{Fixed-mesh free-surface iteration with optional adaptivity}
\label{alg:fs_adapt}
\begin{algorithmic}[1]
\State Initialize $\Omega_w^{(0)}$, set $\Omega_d^{(0)}=\Omega\setminus\Omega_w^{(0)}$, choose $\alpha \ll 1$ and tolerances $\varepsilon_x$
\For{$\ell=0,1,\dots$} \Comment{optional outer adaptivity loop}
    \State Build $\mathbf{k}^{(0)}(\mathbf{x})$ from $\Omega_w^{(0)}$ using Eq.~\eqref{eq:fs_k_update}
    \For{$m=0,1,\dots$} \Comment{inner fixed-mesh iteration}
        \State Impose boundary conditions on $\Gamma_1$, $\Gamma_2$, $\Gamma_5$, and the current seepage-face segment on $\Gamma_4$
        \State Solve Eq.~\eqref{eq:fs_gov} on the fixed mesh using the polygonal CSFEM to obtain $h^{(m)}$
        \State Compute $\psi^{(m)}=h^{(m)}-z$ and update $\Omega_w^{(m+1)}$ and the seepage face using Eq.~\eqref{eq:fs_sat_update}
        \State Update $\mathbf{k}^{(m+1)}(\mathbf{x})$ using Eq.~\eqref{eq:fs_k_update}
        \If{convergence in Eq.~\eqref{eq:fs_stop} is satisfied}
            \State \textbf{break}
        \EndIf
    \EndFor
    \State \textbf{if adaptivity is not used:} \textbf{stop} and return the converged solution
    \State Compute $\eta_K$ using Eqs.~\eqref{eq:band_indicator}--\eqref{eq:eta_indicator} and mark elements by Eq.~\eqref{eq:dorfler_marking}
    \If{$\|\eta\|_{\ell^2(\mathcal{T}_h)}$ in Eq.~\eqref{eq:eta_norm} is below the prescribed tolerance or $\mathcal{M}=\emptyset$}
        \State \textbf{stop} and return the converged solution
    \EndIf
    \State Locally refine the quadtree mesh with 2:1 balance and rebuild hybrid polygonal cells around hanging nodes
    \State Transfer the converged $h$ to the refined mesh by interpolation, set $\Omega_w^{(0)}$ from the transferred $\psi=h-z$, and continue
\EndFor
\end{algorithmic}
\end{algorithm}

The adaptive refinement step is particularly convenient in the present framework because hanging nodes created by local quadtree refinement are naturally handled by the hybrid quadtree discretization without additional constraint equations or transition elements.

\section{Numerical examples}
\label{numerical_examples}
This section presents several benchmark examples to verify the accuracy and convergence behavior of the proposed framework in two-dimensional seepage simulations. To evaluate the performance of the polygonal CSFEM, the obtained results are systematically compared with those from conventional finite element analyses conducted using ABAQUS. All computations were performed on a workstation equipped with an Intel Core i7-4710MQ processor (2.50~GHz) and 4~GB of RAM. The accuracy of the proposed approach is assessed through the relative error of the hydraulic head, defined as:
\begin{equation}
e_{L_2} =
\frac{\left\| \mathbf{H}_{num} - \mathbf{H}_{ref} \right\|_{L_2}}
{\left\| \mathbf{H}_{ref} \right\|_{L_2}},
\end{equation}
where $\mathbf{H}_{num}$ denotes the computed hydraulic head obtained from the polygonal CSFEM, and $\mathbf{H}_{ref}$ corresponds to the reference or exact solution.

\subsection{Patch test}

\begin{figure}[H]
    \centering
    \includegraphics[width=1.0\linewidth]{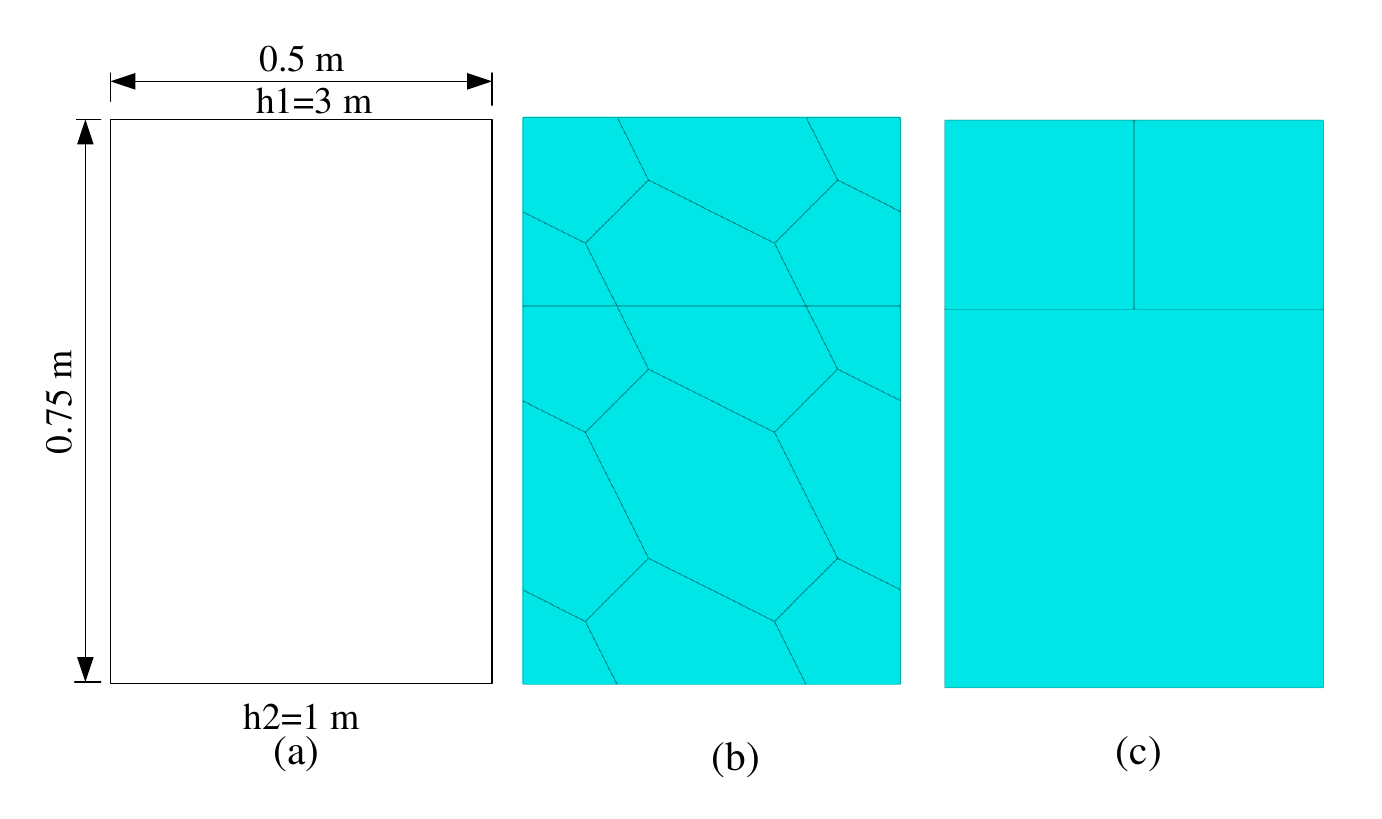}
    \caption{Geometry and mesh models of the patch test: (a) geometry and boundary conditions; (b) polygonal mesh; (c) quadtree mesh.}
    \label{fig:patch_mesh}
\end{figure}

A standard patch test was conducted to verify that the proposed formulation satisfies the fundamental convergence requirement. The geometry and boundary conditions are shown in Fig.~\ref{fig:patch_mesh}(a). Dirichlet boundary conditions with hydraulic heads of 3~m and 1~m were applied at the top and bottom boundaries, respectively. Two mesh configurations were considered, namely the polygonal and quadtree meshes shown in Figs.~\ref{fig:patch_mesh}(b) and (c). The hydraulic conductivity was set to $k = 1 \times 10^{-5}~\mathrm{m/s}$.

\begin{figure}[H]
    \centering
    \includegraphics[width=1\linewidth]{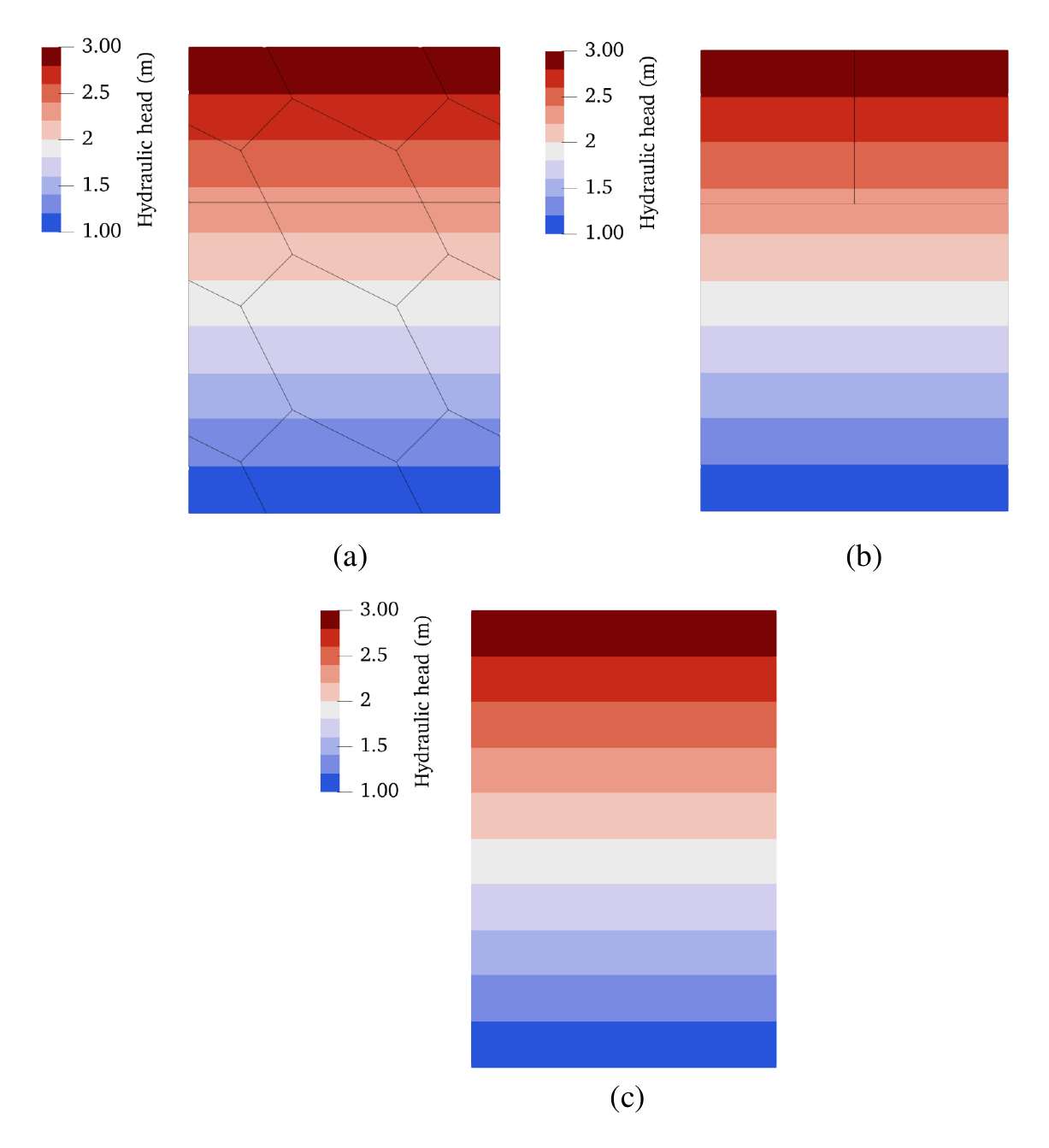}
    \caption{Hydraulic head distribution for the patch test: (a) polygonal mesh; (b) quadtree mesh; (c) exact solution.}
    \label{fig:patch_contour}
\end{figure}

Fig.~\ref{fig:patch_contour} presents the hydraulic head distributions for both mesh types. The contour lines exhibit a perfectly linear and uniform gradient, which agrees with the analytical solution. This indicates that the proposed CSFEM formulation can reproduce linear fields with high accuracy and is insensitive to mesh topology, thereby demonstrating excellent numerical stability and convergence characteristics. Tab.~\ref{tab:patch_Hydraulic_head_value} reports the mesh statistics and the relative errors for the two discretizations. Both meshes achieve errors on the order of $10^{-8}$, confirming that the proposed method successfully passes the patch test. 

\begin{table}[H]
    \centering
    \caption{Comparison of relative error between the polygonal and quadtree meshes.}
    \begin{tabular}{lccc}
        \toprule             
        Mesh type & Elements & Nodes & Relative error \\
        \midrule        
        Polygonal mesh & 15 & 30 & $1.85\times10^{-8}$\\
        Quadtree mesh & 3 & 8 & $1.28\times10^{-8}$\\
        \bottomrule       
    \end{tabular}
    \label{tab:patch_Hydraulic_head_value}
\end{table}

\subsection{Steady-state seepage analysis}
\subsubsection{Foundation of a concrete dam}
\begin{figure}[H]
    \centering
    \includegraphics[width=0.8\linewidth]{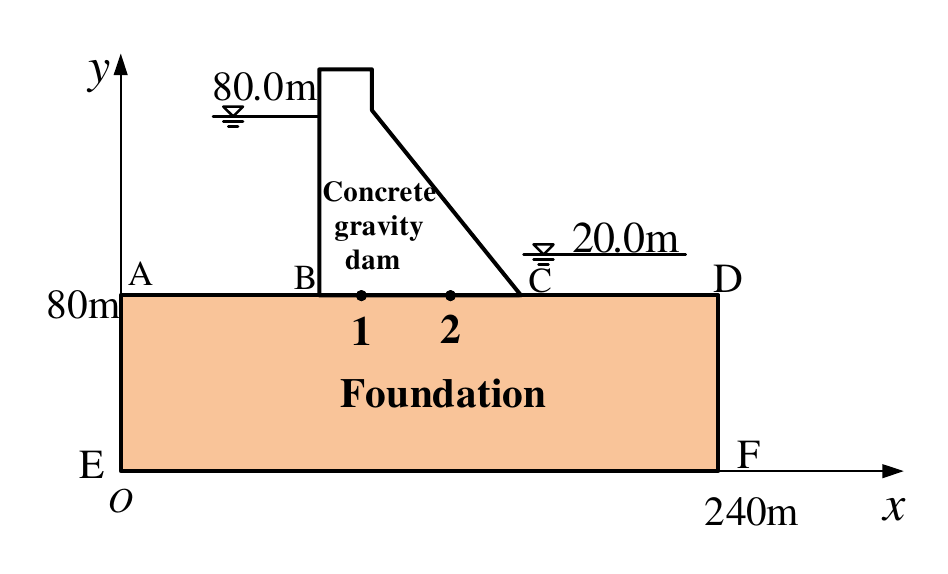}
    \caption{Geometry and boundary conditions of the concrete dam.}
    \label{fig:ex01_geo}
\end{figure}

In this example, a steady-state seepage problem of a dam foundation is analyzed, as illustrated in Fig.~\ref{fig:ex01_geo}. The dam is assumed to be impervious, with boundaries BC, DE, FE, and EA defined as impermeable surfaces. A hydraulic head of $80~\mathrm{m}$ is applied along boundary AB, while boundary CD is subjected to a hydraulic head of $20~\mathrm{m}$. To verify the accuracy of the proposed method, two monitoring points are selected: point~1~$(100,\,80)$ and point~2~$(140,\,80)$. The permeability coefficient of the dam foundation is $k_x = k_y = 1 \times 10^{-5}~\mathrm{cm/s}$. For comparison, the commercial software ABAQUS employs CPE4P elements, whereas the proposed CSFEM utilizes polygonal elements.

\begin{figure}[H]
    \centering
    \includegraphics[width=0.8\linewidth]{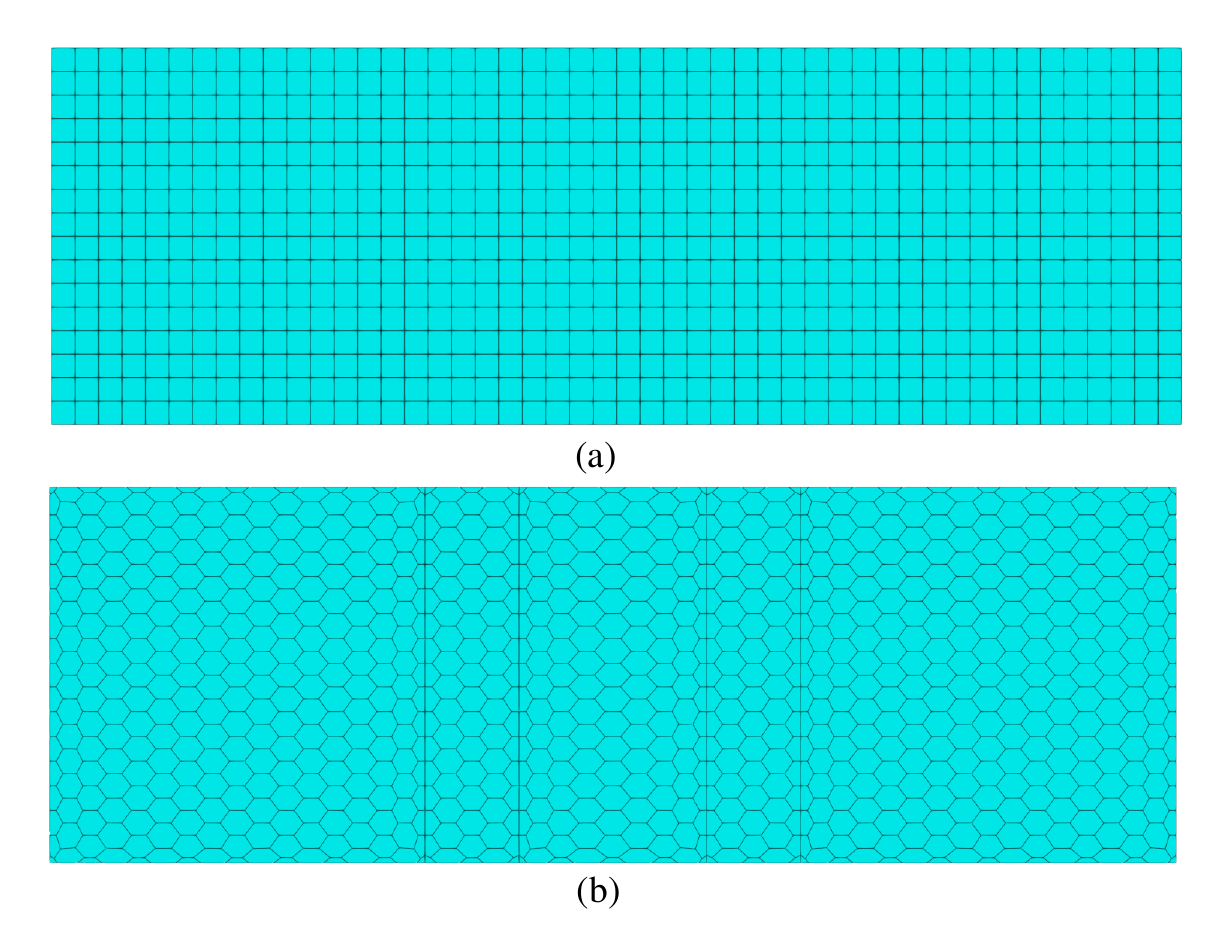}
    \caption{Mesh models of the dam foundation: (a) FEM mesh; (b) CSFEM mesh.}
    \label{fig:ex01_mesh}
\end{figure}

\begin{figure}[H]
    \centering
    \begin{subfigure}[b]{0.48\linewidth}
        \centering
        \includegraphics[width=\linewidth]{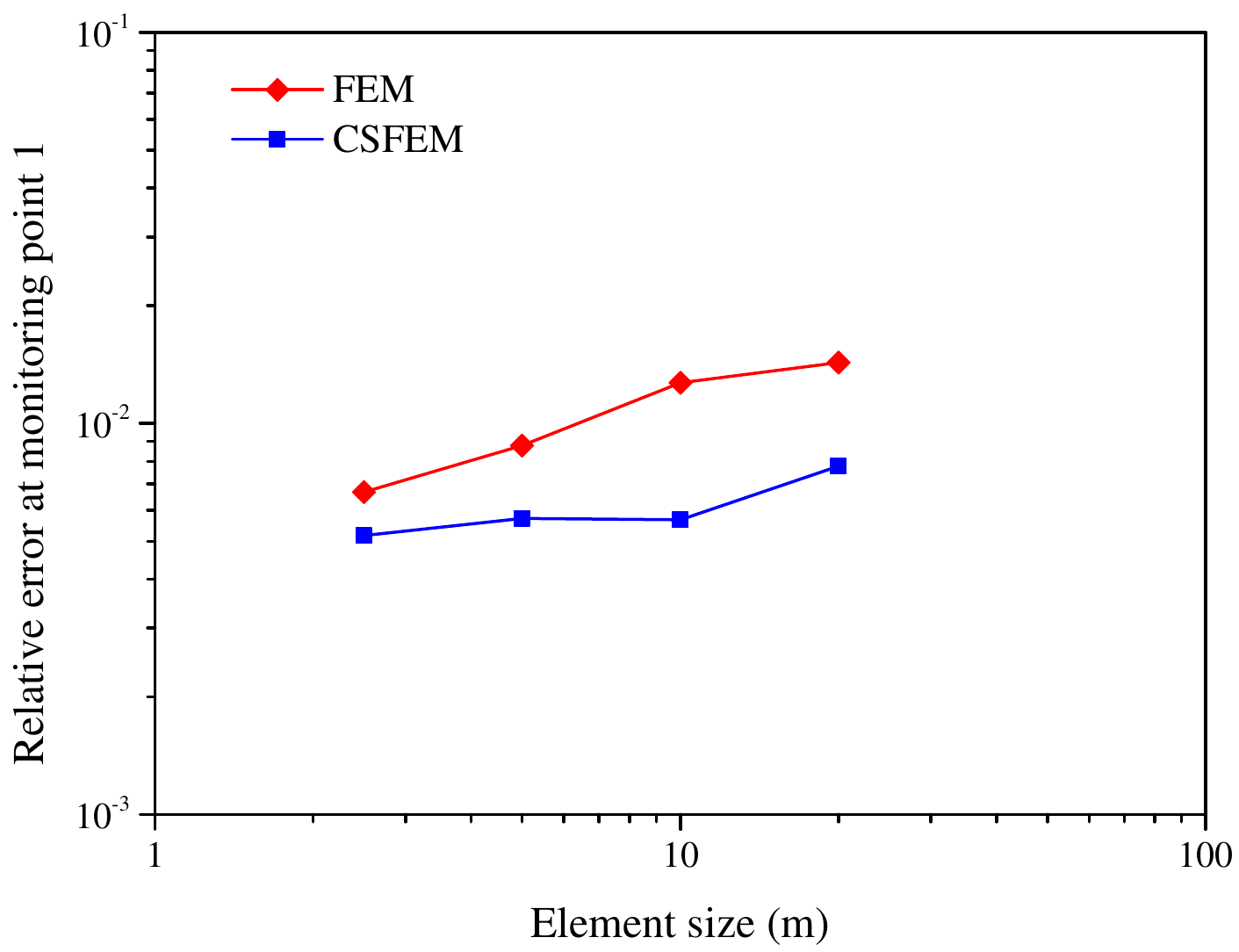}
        \caption{}
    \end{subfigure}
    \hfill
    \begin{subfigure}[b]{0.48\linewidth}
        \centering
        \includegraphics[width=\linewidth]{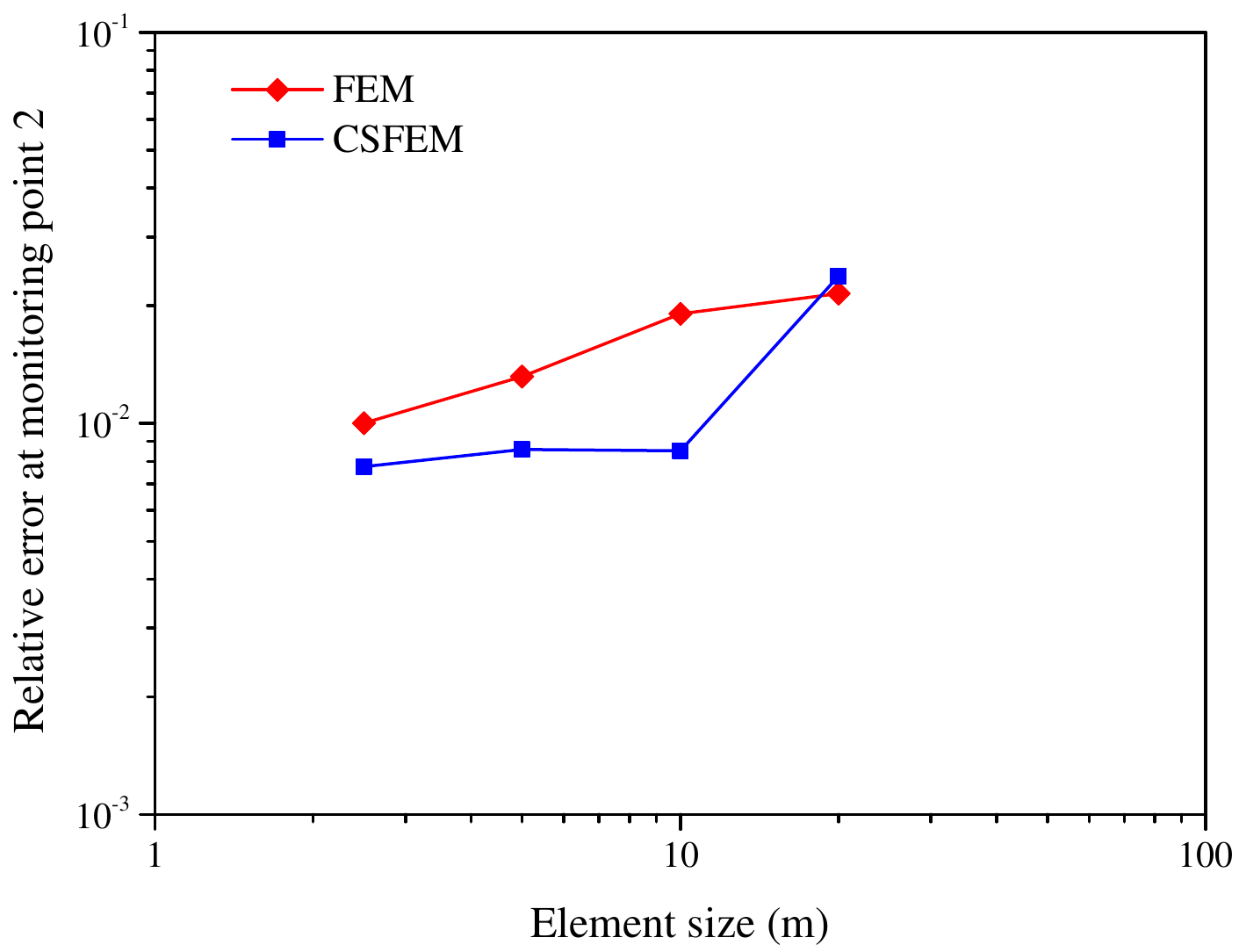}
        \caption{}
    \end{subfigure}
    \caption{Comparison of convergence rates of hydraulic head at the monitoring points: (a) point 1; (b) point 2.}
    \label{fig:ex01_error}
\end{figure}

\begin{figure}[H]
    \centering
    \includegraphics[width=0.8\linewidth]{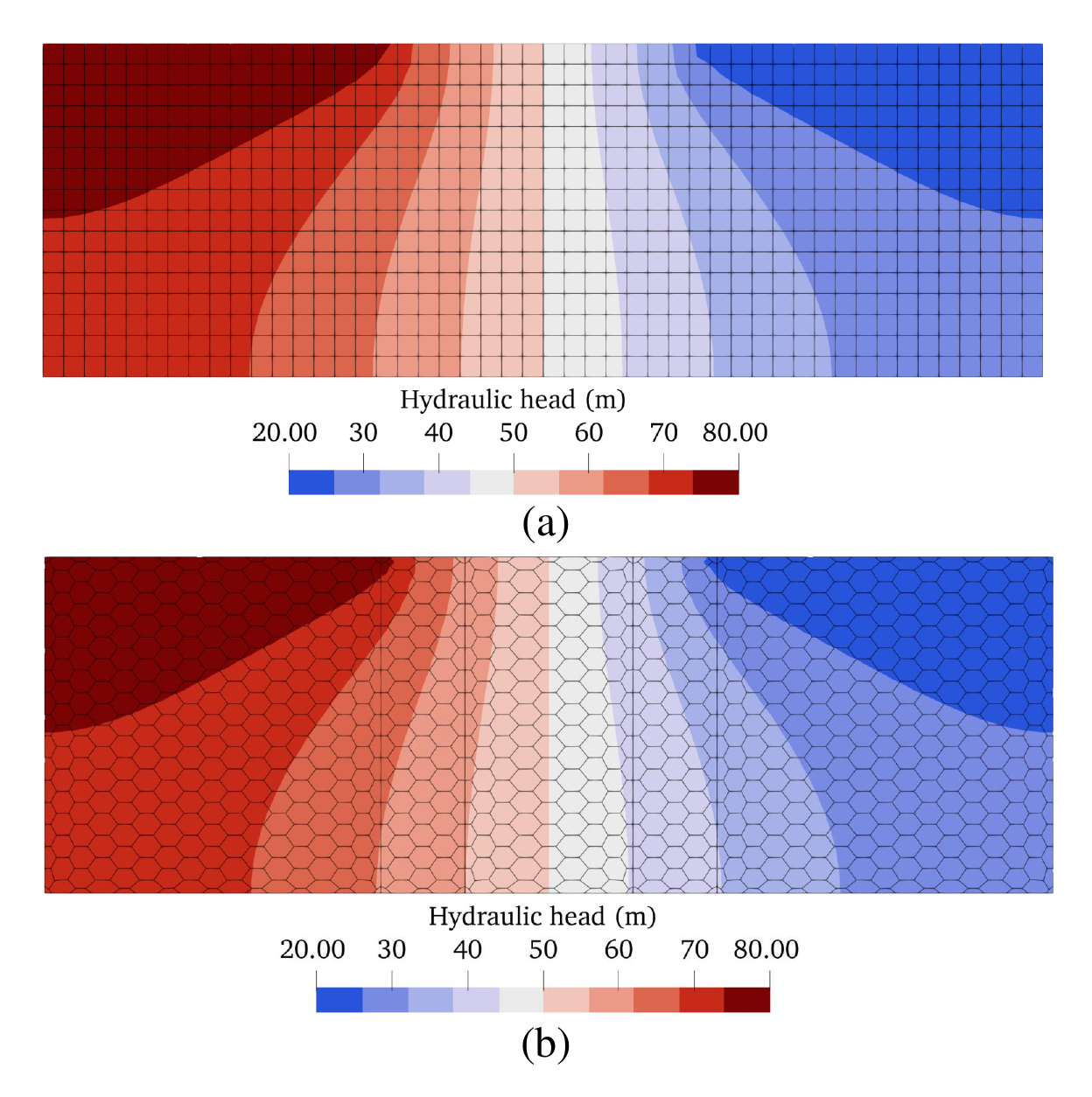}
    \caption{Hydraulic head distribution in the dam foundation: (a) FEM; (b) CSFEM.}
    \label{fig:ex01_contour}
\end{figure}

A convergence study is performed through mesh $h$-refinement with mesh sizes of $20~\mathrm{m}$, $10~\mathrm{m}$, $5~\mathrm{m}$, and $2.5~\mathrm{m}$. The mesh with an element size of $5~\mathrm{m}$ is shown in Fig.~\ref{fig:ex01_mesh}. Fig.~\ref{fig:ex01_error} presents the relative errors of hydraulic head at the monitoring points, indicating that the polygonal CSFEM and the conventional FEM exhibit comparable convergence rates under mesh refinement. Tab.~\ref{tab:ex01_t1} lists the hydraulic head values for the $5~\mathrm{m}$ mesh. The relative errors for the CSFEM and FEM are $7.1\times10^{-3}$ and $1.1\times10^{-2}$, respectively, indicating that the CSFEM achieves higher computational accuracy. Furthermore, Fig.~\ref{fig:ex01_contour} shows that the hydraulic head contours obtained by the CSFEM are in excellent agreement with those of the FEM.

\begin{table}[H]
    \centering
    \small
    \caption{Comparison of hydraulic head at monitoring points obtained by the CSFEM and FEM (element size: 5 m).}
    \label{tab:ex01_t1}
    \begin{tabular}{cccc}
        \hline
        Method & Point 1 (m) & Point 2 (m) & Relative error \\ 
        \hline
        Analytical solution \citep{Li2003Comparisons} & 60.00 & 40.00 & -- \\
        CSFEM & 60.3428 & 39.6572 & $7.1\times10^{-3}$ \\
        FEM  & 60.5265 & 39.4735 & $1.1\times10^{-2}$ \\
        \hline
    \end{tabular}
\end{table}

\subsubsection{Permeable media}
\begin{figure}[H]
    \centering
    \includegraphics[width=0.8\linewidth]{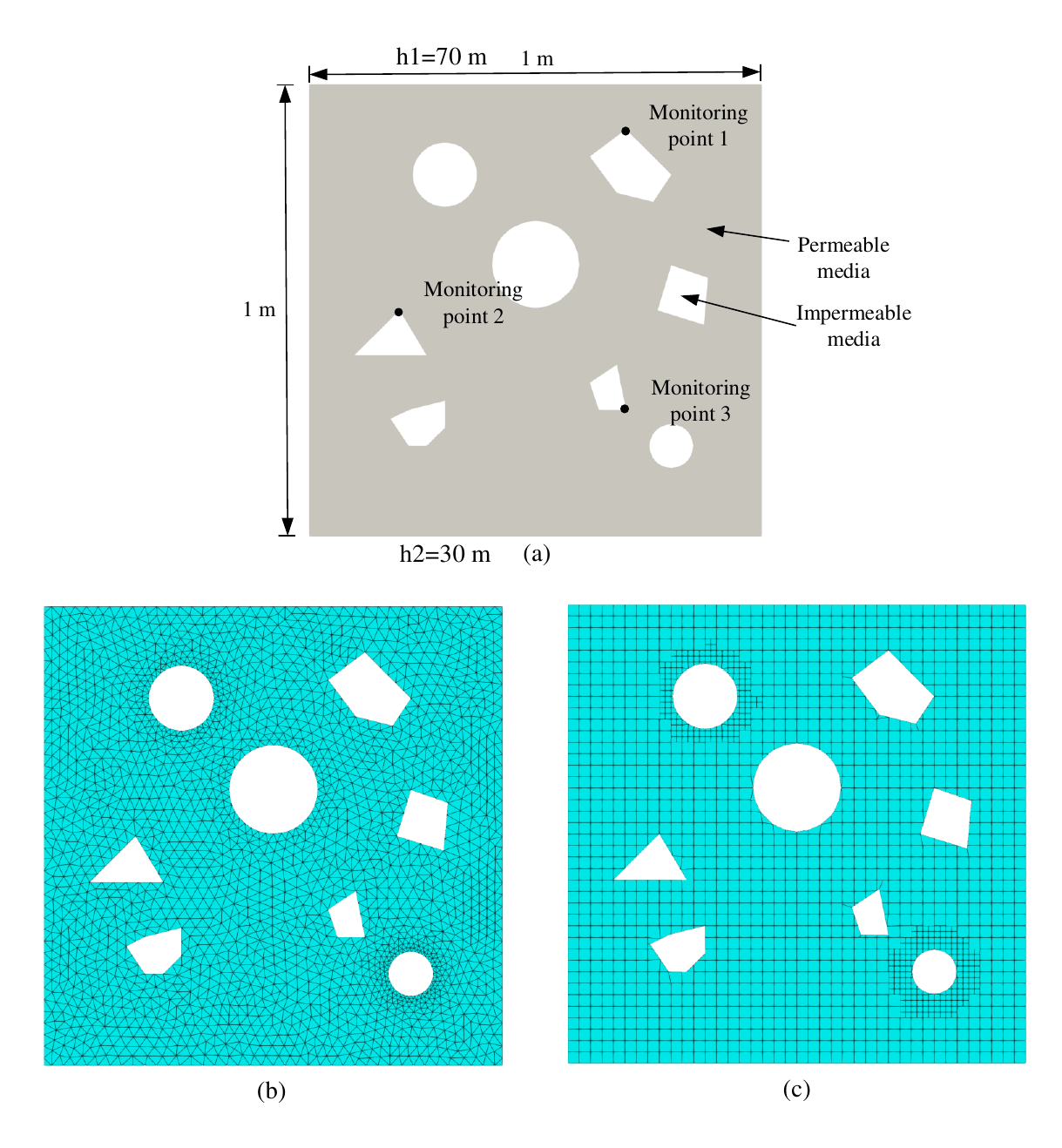}
    \caption{Geometric and mesh model of the permeable medium: (a) geometric model and boundary conditions; (b) FEM mesh; (c) hybrid mesh (CSFEM).}
    \label{fig:ex02_geo_mesh}
\end{figure}

To evaluate the flexibility of the proposed CSFEM combined with a hybrid quadtree mesh, a steady-state seepage problem is analyzed for a permeable medium containing an impermeable inclusion, as illustrated in Fig.~\ref{fig:ex02_geo_mesh}(a). The domain measures 1 m in both length and width, and the permeability coefficient is $k = 2\times10^{-6}$ cm/s. A hydraulic head of $h_1 = 70$ m is applied along the top boundary, and $h_2 = 30$ m along the bottom boundary.

The hybrid quadtree mesh, shown in Fig.~\ref{fig:ex02_geo_mesh}(c), allows adaptive refinement near curved boundaries, thereby improving boundary conformity without increasing the overall mesh density. Regular regions are discretized using standard quadtree elements, whereas irregular hybrid elements are employed along curved boundaries to achieve accurate geometric representation. For comparison, the conventional FEM adopts an unstructured triangular mesh to model the same geometry, as shown in Fig.~\ref{fig:ex02_geo_mesh}(b). The reference solution is obtained from a FEM analysis using a sufficiently refined mesh.

\begin{figure}[H]
    \centering
    \includegraphics[width=0.8\linewidth]{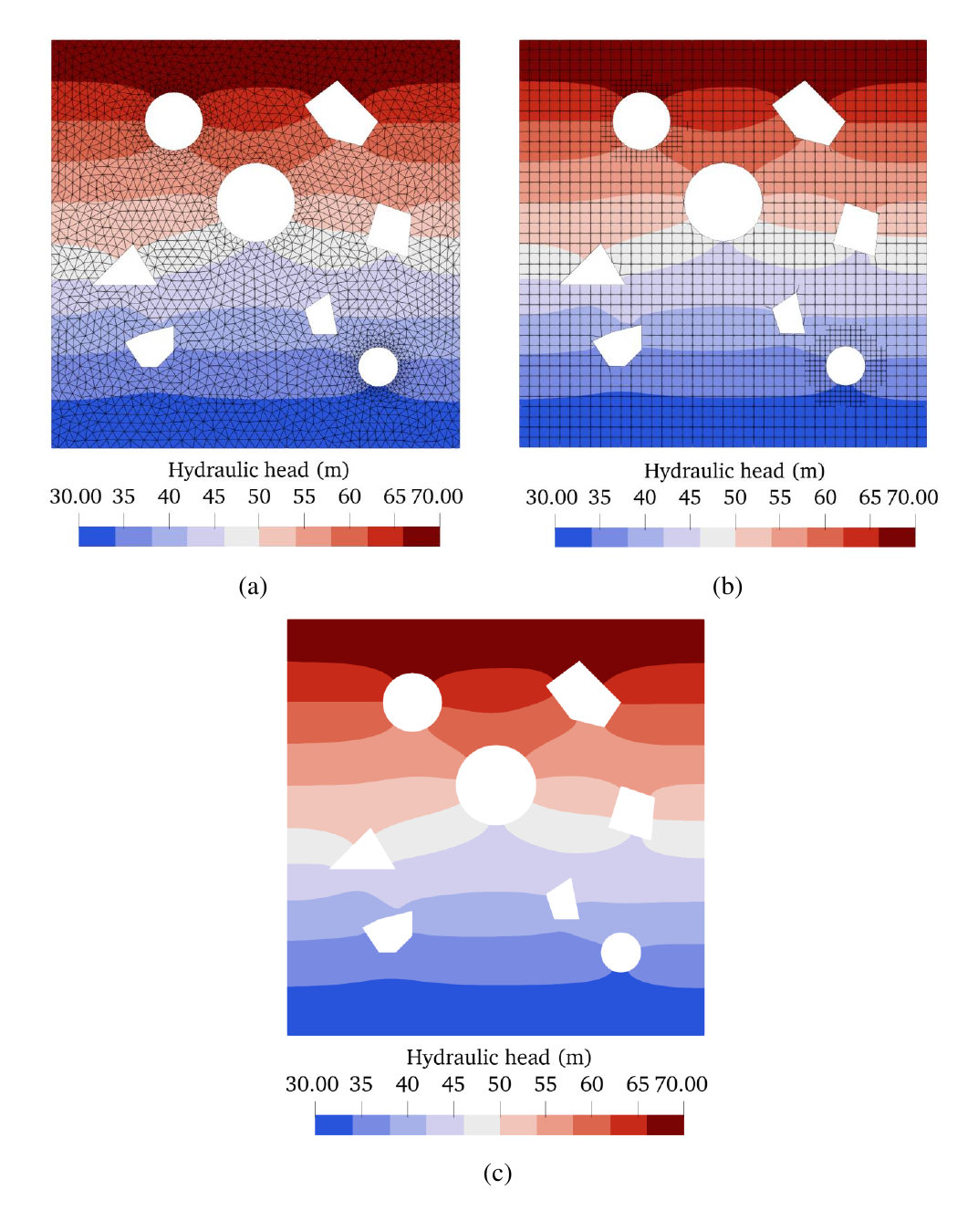}
    \caption{Hydraulic head distribution in the permeable medium: (a) FEM; (b) CSFEM; (c) reference solution.}
    \label{fig:ex02_contour}
\end{figure}

To assess the computational accuracy, three monitoring points are placed within the domain, as shown in Fig.~\ref{fig:ex02_geo_mesh}(c). Tab.~\ref{tab:ex02_t1} summarizes the hydraulic head results at these points obtained by the CSFEM and FEM. The error of the CSFEM is $5.01\times10^{-4}$, while that of the FEM is $6.05\times10^{-4}$, demonstrating that the CSFEM achieves slightly higher accuracy. Furthermore, Fig.~\ref{fig:ex02_contour} presents the hydraulic head distribution, showing that the CSFEM solution is in excellent agreement with both the FEM and the reference solution.

\begin{table}[H]
    \centering
    \small
    \caption{Comparison of hydraulic head at three monitoring points obtained by the CSFEM and FEM.}
    \label{tab:ex02_t1}
    \begin{tabular}{ccccc}
        \hline
        Method & Point 1 (m) & Point 2 (m) & Point 3 (m) & Relative error \\ 
        \hline
        Reference & 67.9313 & 50.7946 & 39.7161 & -- \\
        CSFEM & 67.9679 & 50.8353 & 39.7226 & $5.01\times 10^{-4}$ \\
        FEM & 67.8952 & 50.7381 & 39.7229 & $6.05\times 10^{-4}$ \\
        \hline
    \end{tabular}
\end{table}

\subsection{Transient seepage analysis}
\subsubsection{Dam foundation with irregular geometry}
\begin{figure}[H]
    \centering
    \includegraphics[width=1.0\linewidth]{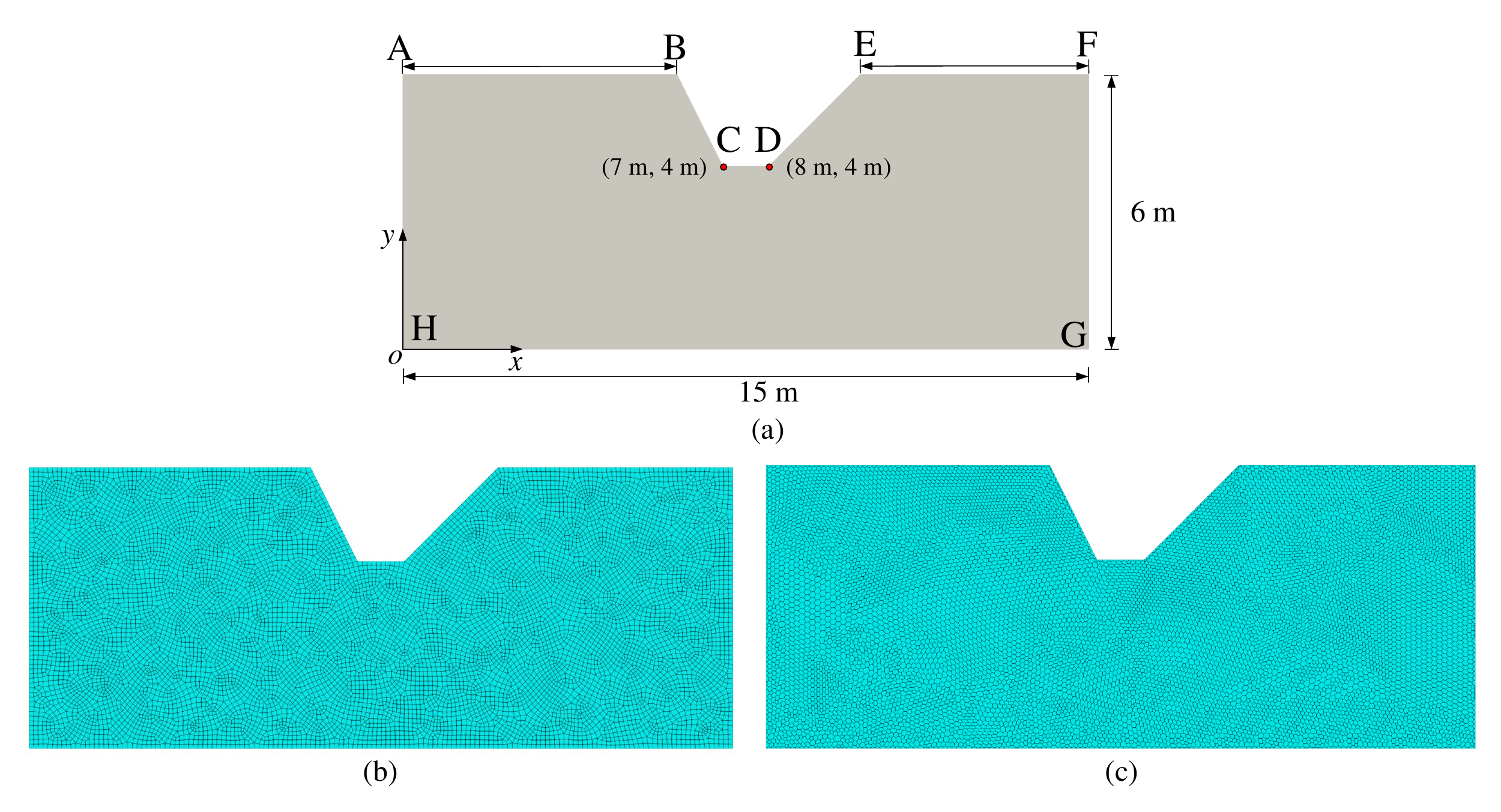}
    \caption{Geometric and mesh model of dam foundation: (a) geometric model and boundary conditions; (b) FEM mesh; (c) polygonal mesh (CSFEM).}
    \label{fig:ex03_geo_mesh}
\end{figure}

\begin{figure}[H]
    \centering
    \includegraphics[width=0.8\linewidth]{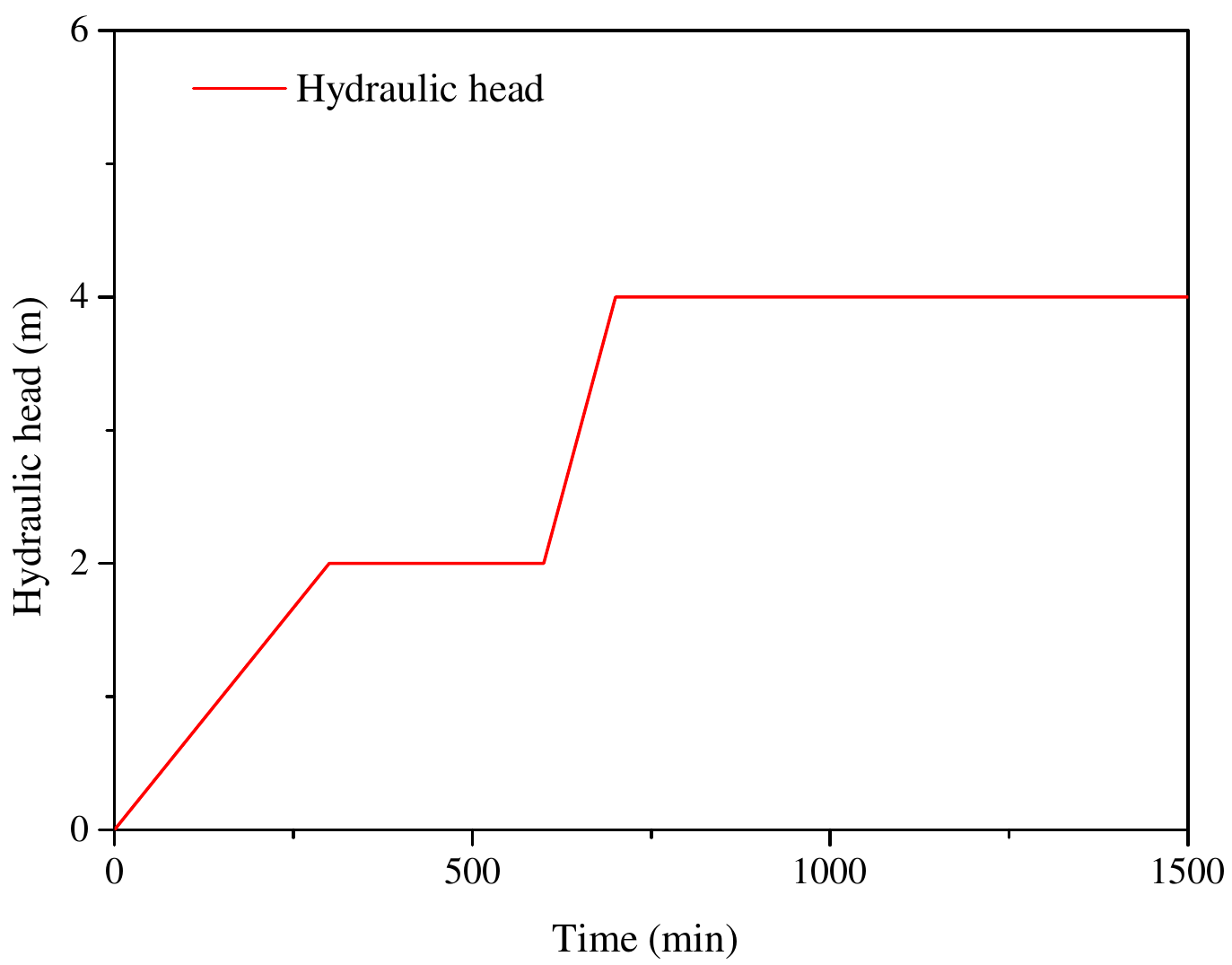}
    \caption{Hydraulic head boundary condition in the upstream.}
    \label{fig:ex03_hyheadhis}
\end{figure}

In this example, a concrete dam with an irregular foundation is analyzed, as illustrated in Fig.~\ref{fig:ex03_geo_mesh}. To investigate the transient seepage behavior, two monitoring locations, C and D, are chosen within the foundation domain. The storage coefficient is taken as $S_s = 1.0 \times 10^{-3}~\mathrm{m^{-1}}$, and the hydraulic conductivity is $k = 1.0 \times 10^{-3}~\mathrm{m/min}$. The initial upstream and downstream water levels are 0.0 m and 0.0 m, respectively. Over time, the upstream hydraulic head is progressively increased from 0.0 m to 2.0 m, and eventually reaches 4.0 m, as illustrated in Fig.~\ref{fig:ex03_hyheadhis}. A reference solution is obtained using a sufficiently refined FEM mesh.

\begin{figure}[H]
    \centering
    \begin{subfigure}[b]{0.48\linewidth}
        \centering
        \includegraphics[width=\linewidth]{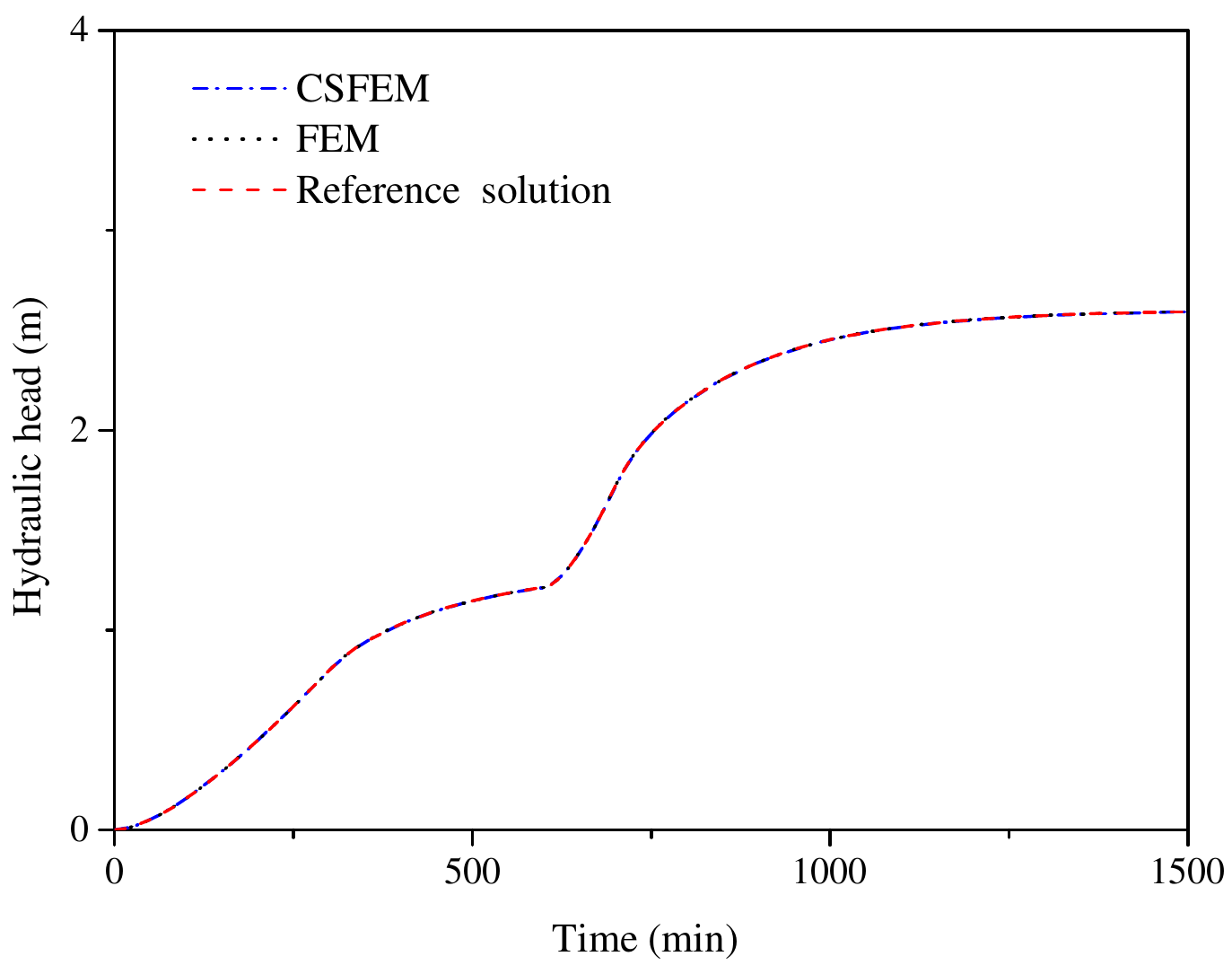}
        \caption{}
    \end{subfigure}
    \hfill
    \begin{subfigure}[b]{0.48\linewidth}
        \centering
        \includegraphics[width=\linewidth]{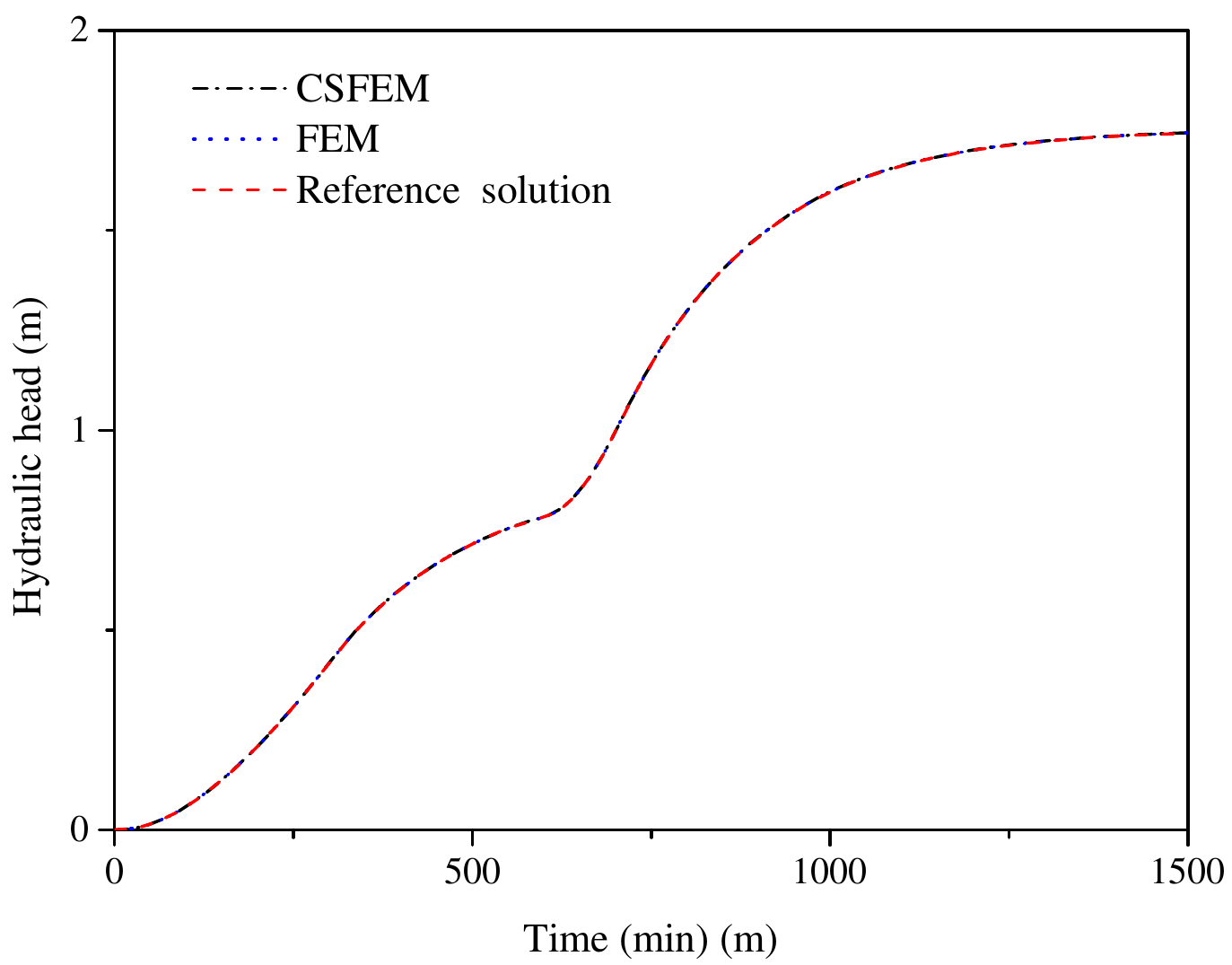}
        \caption{}
    \end{subfigure}
    \caption{Temporal evolution of the hydraulic head at monitoring points: (a) point C; (b) point D.}
    \label{fig:ex03_monitoring_his}
\end{figure}

\begin{figure}[H]
    \centering
    \includegraphics[width=1.0\linewidth]{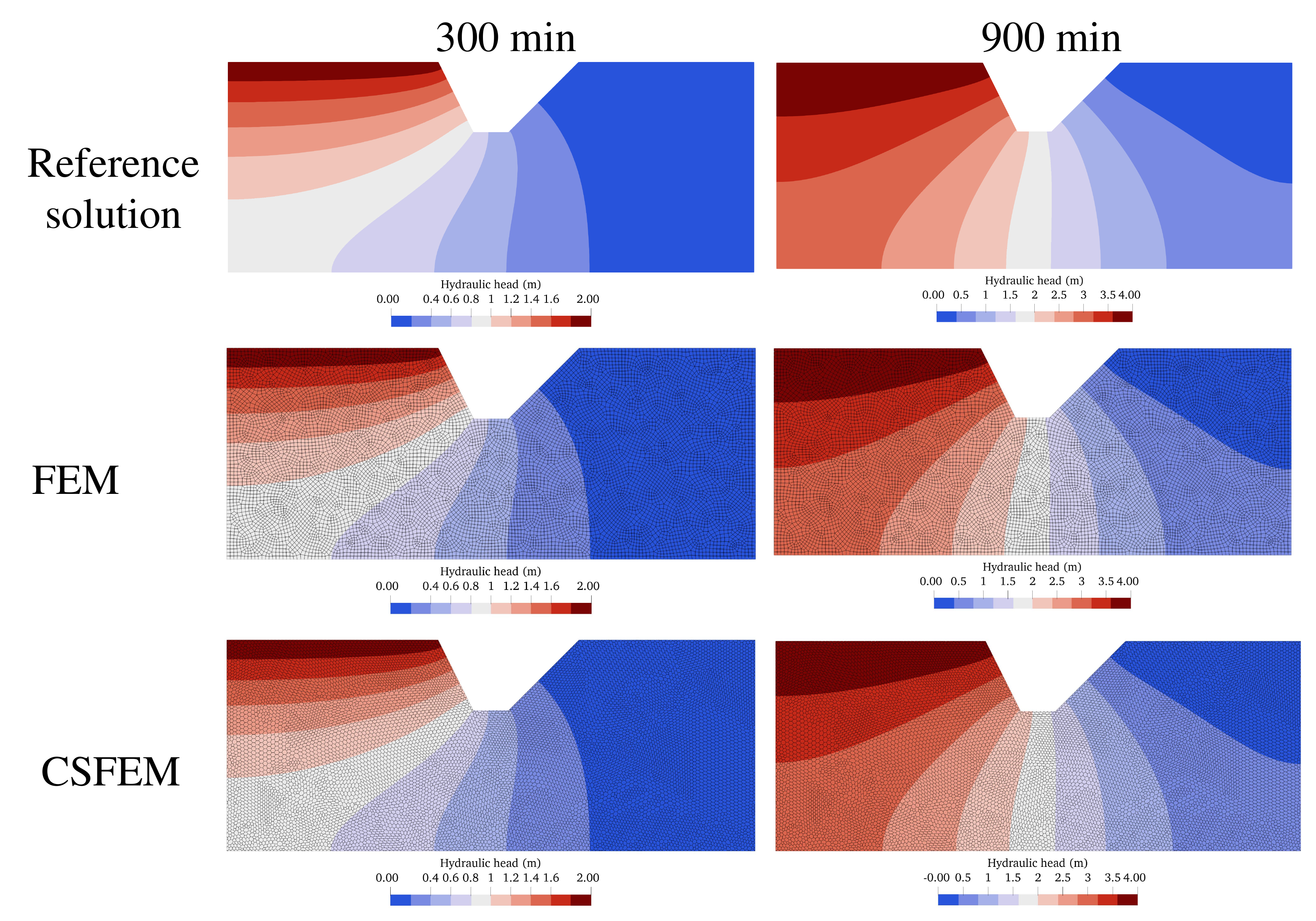}
    \caption{Hydraulic head distribution in the dam foundation.}
    \label{fig:ex03_contour}
\end{figure}

Fig.~\ref{fig:ex03_monitoring_his} shows the transient hydraulic head responses at monitoring points C and D, where the CSFEM predictions are in close agreement with the reference solution and the conventional FEM, capturing the gradual head rise induced by the staged upstream boundary condition. Fig.~\ref{fig:ex03_contour} presents the hydraulic head contours at the final time, demonstrating a smooth distribution consistent with the FEM reference and confirming the robustness of the proposed CSFEM on irregular dam-foundation geometry.

\subsubsection{Transient seepage analysis for complex geometry}
\begin{figure}[H]
    \centering
    \includegraphics[width=1.0\linewidth]{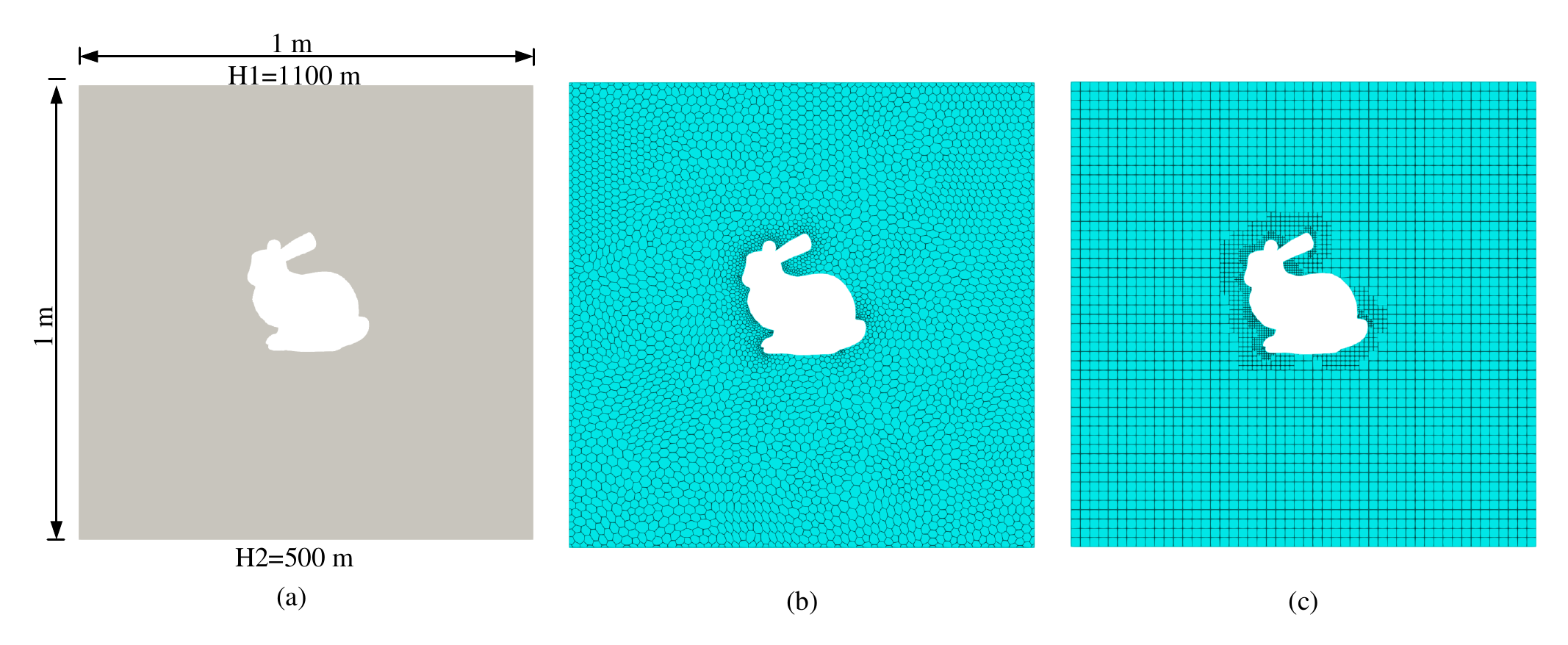}
    \caption{Geometric and mesh model for a square plate with a Stanford bunny cavity; (a) geometric model; (b) polygonal mesh; (c) hybrid quadtree mesh.}
    \label{fig:ex04_geo_mesh}
\end{figure}

A transient seepage problem is considered in a square domain with a side length of $L = 1.0$ m containing a Stanford bunny-shaped cavity~\cite{StanfordBunny1994}, as illustrated in Fig.~\ref{fig:ex04_geo_mesh}(a). The permeability coefficient is set to $k = 2 \times 10^{-7}$ m/s, and the specific storage coefficient is $S_s = 0.001$ m$^{-1}$. The total simulation time is 2000 s with a uniform time increment of $\Delta t = 10$ s for both the CSFEM and FEM analyses. The hydraulic heads at the top and bottom boundaries are prescribed as 1100 m and 500 m, respectively. To assess the performance of the two numerical methods, two monitoring points, denoted as A and B, are selected to compare the transient hydraulic responses.

\begin{figure}[H]
    \centering
    \begin{subfigure}[b]{0.48\linewidth}
        \centering
        \includegraphics[width=\linewidth]{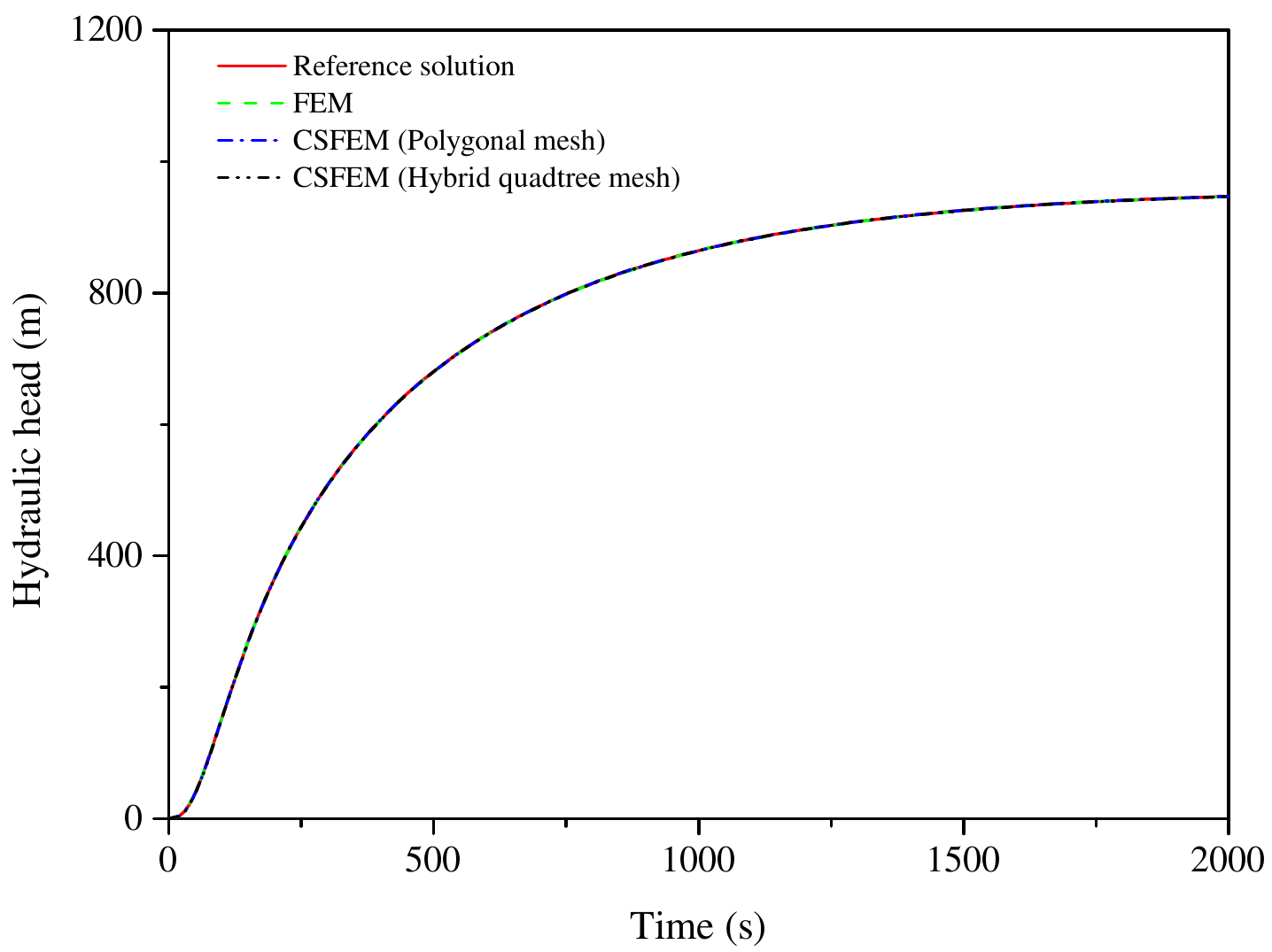}
        \caption{}
    \end{subfigure}
    \hfill
    \begin{subfigure}[b]{0.48\linewidth}
        \centering
        \includegraphics[width=\linewidth]{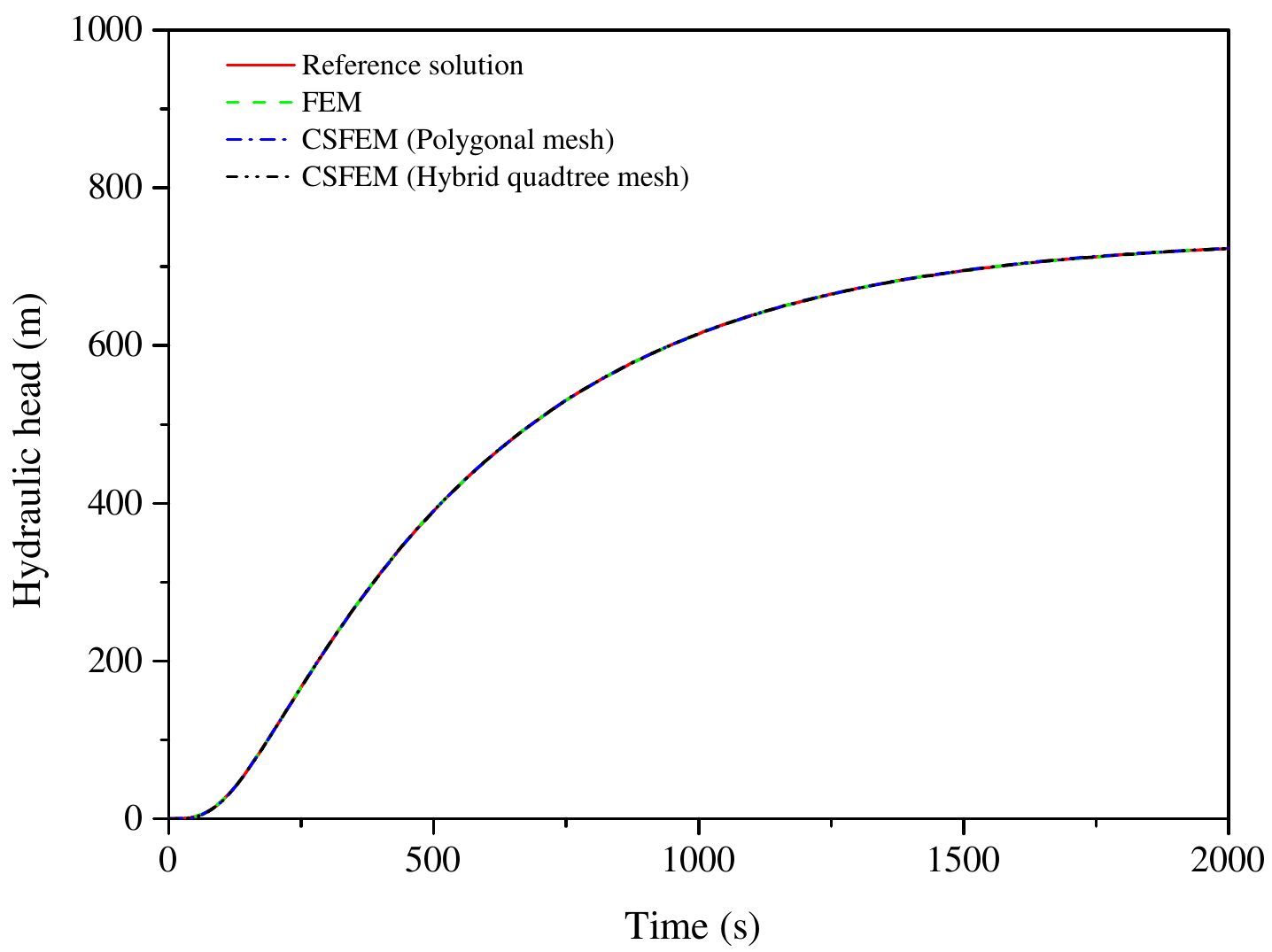}
        \caption{}
    \end{subfigure}
    \caption{Temporal evolution of the hydraulic head at monitoring points: (a) point A; (b) point B.}
    \label{fig:ex04_monitoring_his}
\end{figure}

\begin{figure}[H]
    \centering
    \includegraphics[width=1.0\linewidth]{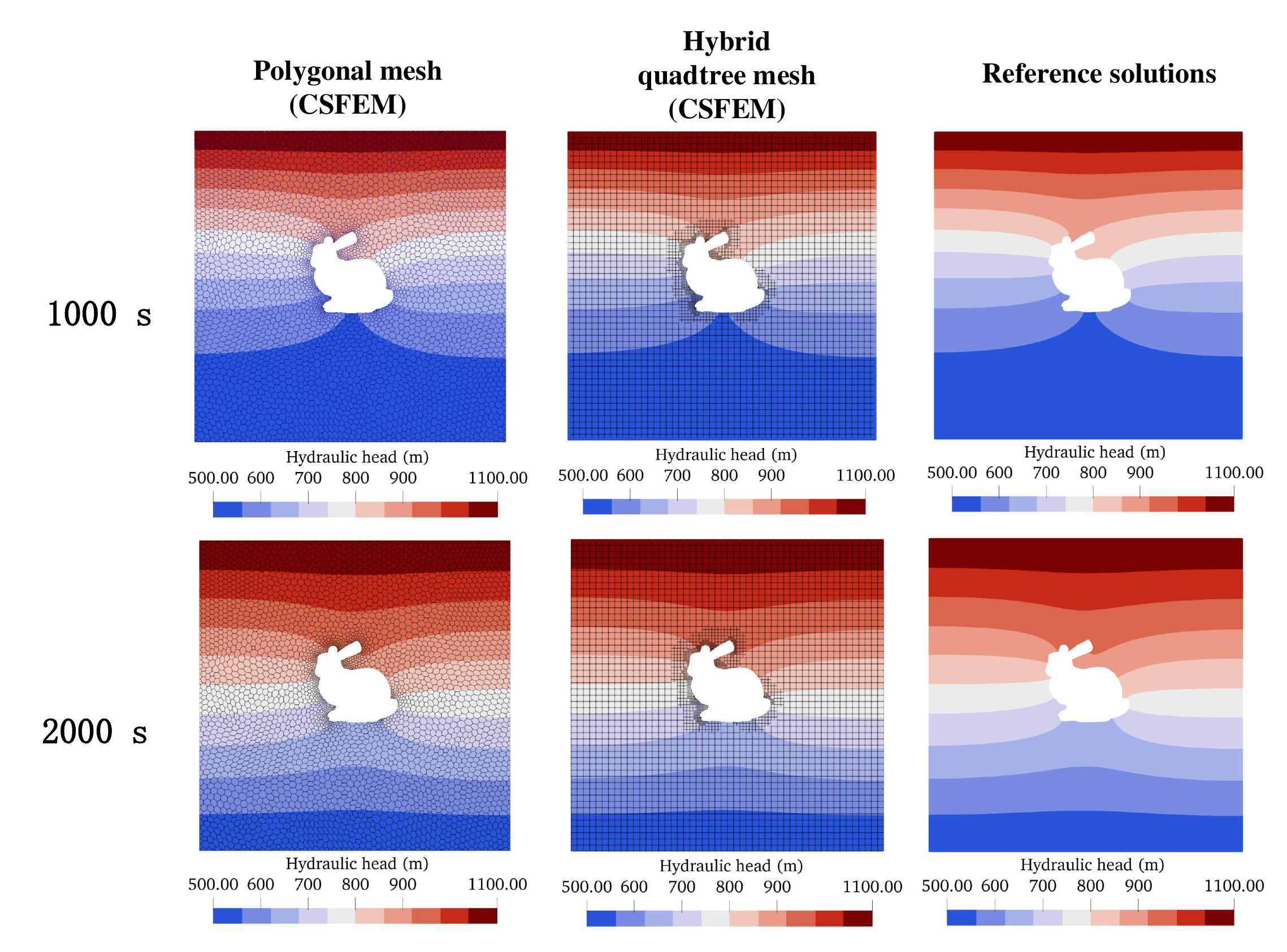}
   \caption{Hydraulic head distributions at different times.}
    \label{fig:ex04_contour}
\end{figure}

Fig.~\ref{fig:ex04_monitoring_his} presents the temporal evolution of the hydraulic head at the monitoring points. The results obtained using the CSFEM with polygonal elements and hybrid quadtree elements exhibit excellent agreement with the reference solution. Moreover, Fig.~\ref{fig:ex04_contour} shows the hydraulic head distribution at different times. Similarly, the hydraulic head distributions obtained by the CSFEM with polygonal elements and hybrid quadtree elements match the reference solution very well.

\subsection{Free-surface seepage problems}
\subsubsection{Free-surface seepage in a homogeneous dam}

\begin{figure}[H]
    \centering
    \includegraphics[width=0.6\linewidth]{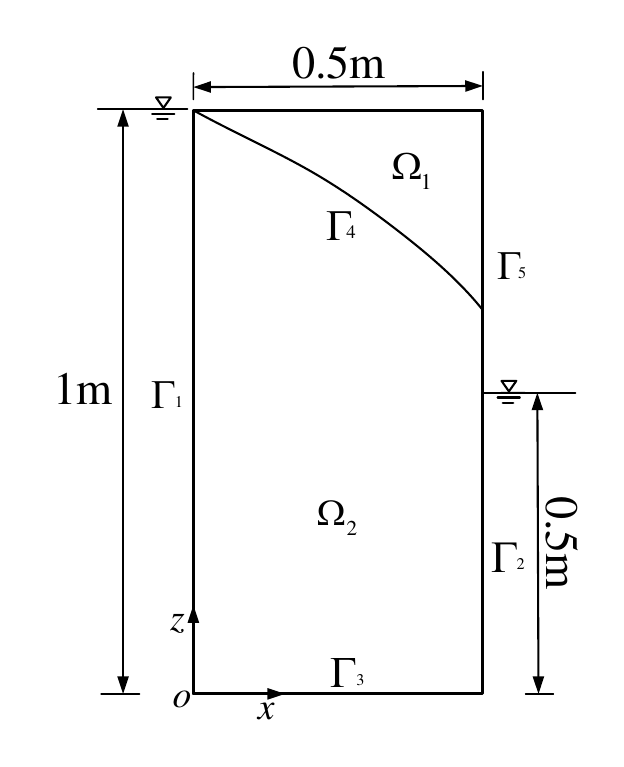}
    \caption{Geometric model of a homogeneous dam.}
    \label{fig:ex05_geo}
\end{figure}

\begin{figure}[H]
    \centering
    \includegraphics[width=0.8\linewidth]{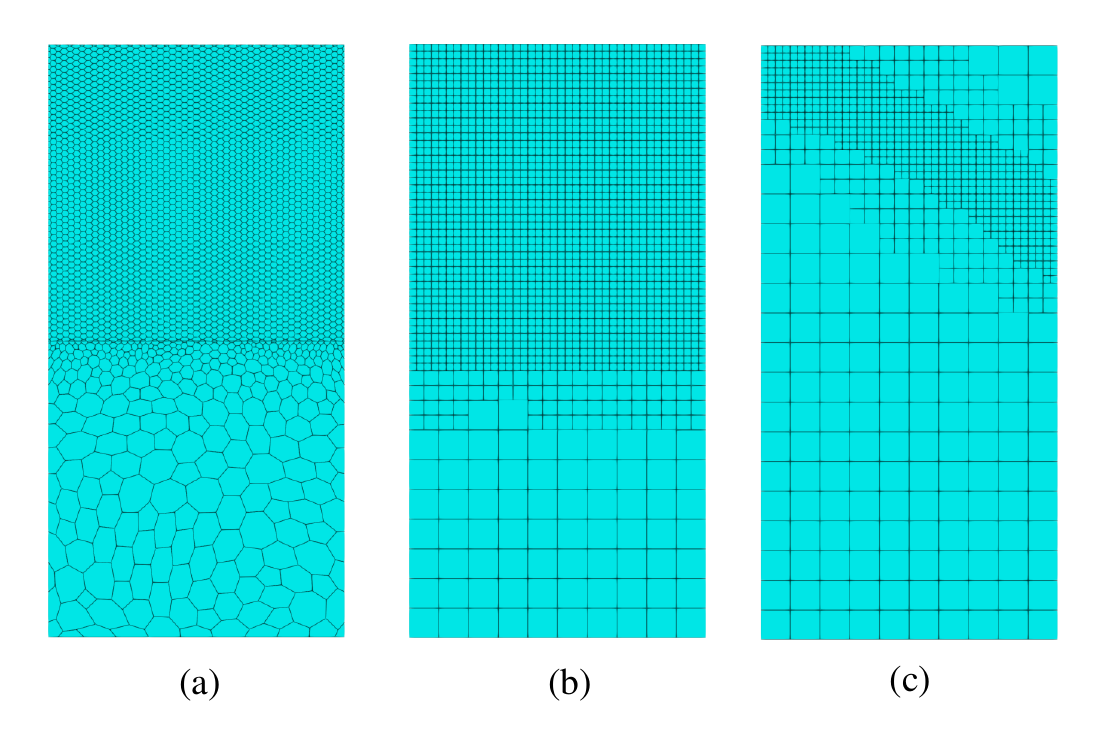}
    \caption{Mesh model of a homogeneous dam; (a) the polygonal mesh; (b) the quadtree mesh; (c) the adaptive mesh. The adaptive mesh in (c) is the final mesh obtained by the indicator in Eqs.~\eqref{eq:band_indicator}--\eqref{eq:grad_indicator} with bulk-chasing marking (Eq.~\eqref{eq:dorfler_marking}, $\theta=0.5$) and 2:1 balanced quadtree refinement.}
    \label{fig:ex05_mesh}
\end{figure}

The benchmark considers a homogeneous rectangular dam, whose geometry is shown in Fig.~\ref{fig:ex05_geo}. The dam has a height of 1.0~m and a base width of 0.5~m. A constant hydraulic head of 1.0~m is prescribed along the upstream boundary, whereas the downstream boundary is maintained at a head of 0.5~m. The foundation boundary is assumed to be impermeable. To adequately capture the free-surface zone forming above the downstream water level, two a priori locally refined meshes with the same refinement pattern are employed, as shown in Fig.~\ref{fig:ex05_mesh}(a) and (b). Specifically, the potential phreatic-zone region is pre-refined to a cell size of 0.0125~m, whereas the remainder of the domain adopts a coarser size of 0.05~m. The conforming polygonal mesh is shown in Fig.~\ref{fig:ex05_mesh}(a), while the corresponding quadtree mesh is shown in Fig.~\ref{fig:ex05_mesh}(b). Moreover, an adaptive quadtree mesh is constructed via solution-driven local refinement in the vicinity of the predicted free-surface region, as illustrated in Fig.~\ref{fig:ex05_mesh}(c), so that resolution is concentrated where steep hydraulic gradients and active-set switching occur while keeping the total degrees of freedom moderate. 

\begin{figure}[H]
    \centering
    \includegraphics[width=0.4\linewidth]{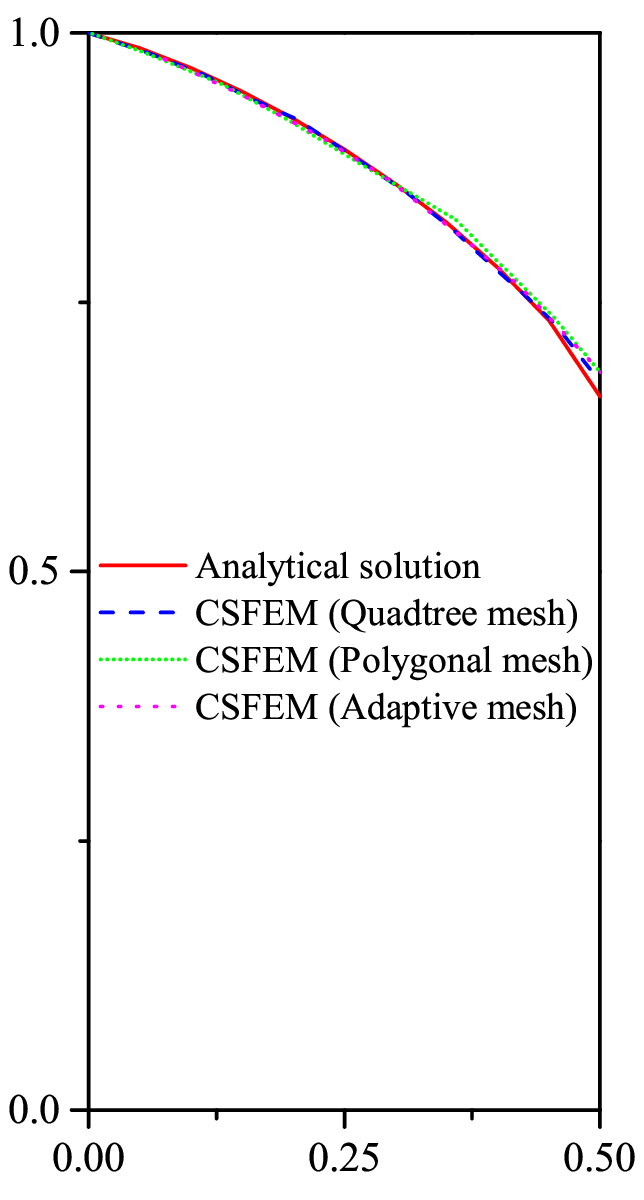}
    \caption{Comparison of free-surface positions of a homogeneous dam.}
    \label{fig:ex05_free_surface}
\end{figure}

Fig.~\ref{fig:ex05_free_surface} presents the free-surface profiles obtained using the proposed CSFEM, together with those from alternative numerical formulations and the analytical reference solution. The three CSFEM results show consistently good agreement with the analytical curve, demonstrating that the method can accurately capture seepage fields involving a moving free surface. In addition, the free-surface curves computed on the quadtree mesh and the adaptive mesh are slightly smoother than that obtained on the polygonal mesh, which can be attributed to the more regular element layout and the locally enhanced resolution achieved by adaptive refinement.

\begin{figure}[H]
    \centering
    \includegraphics[width=1.0\linewidth]{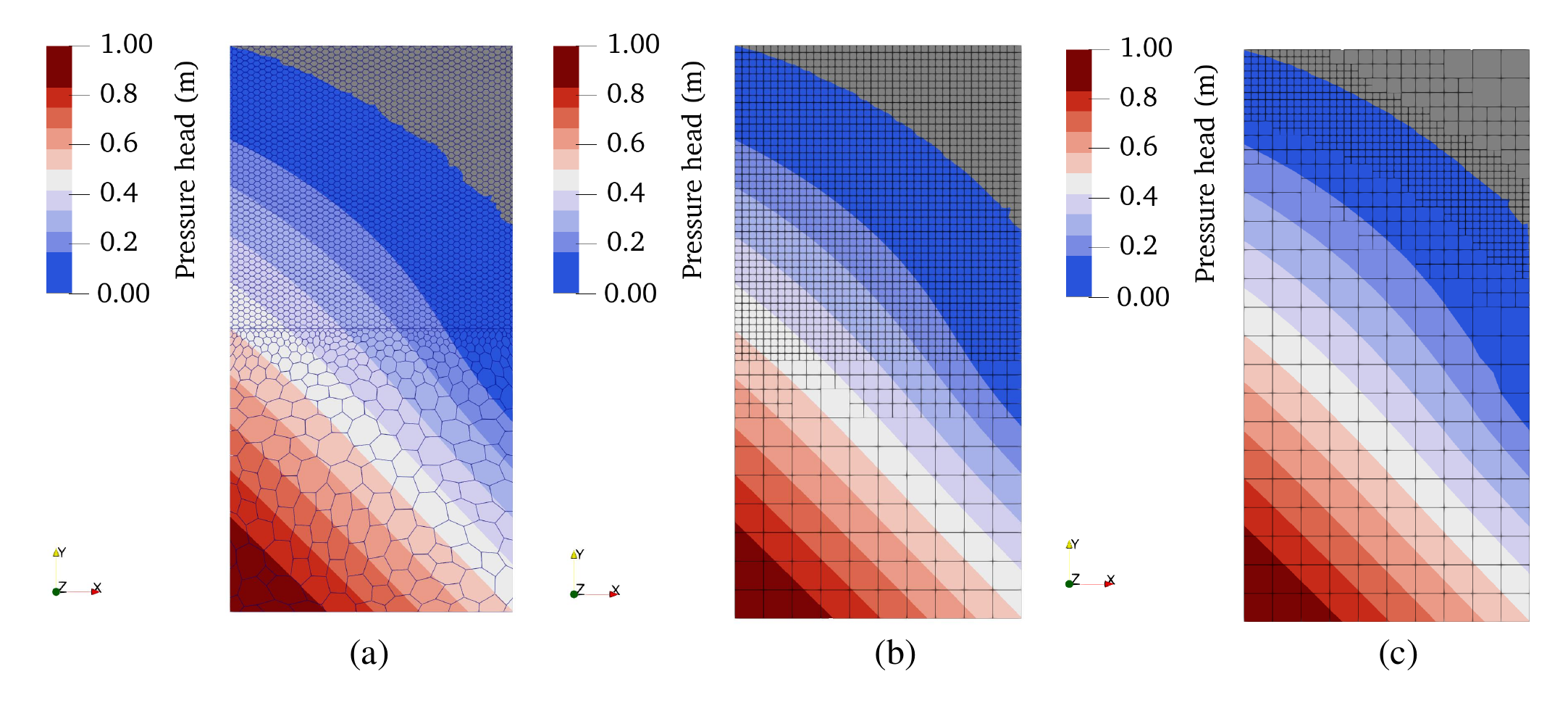}
    \caption{The pressure head distribution within a homogeneous dam; (a) the polygonal mesh; (b) the quadtree mesh; (c) the adaptive mesh.}
    \label{fig:ex05_Pressure_head}
\end{figure}

The overflow-point coordinates obtained using the proposed CSFEM are reported in Tab.~\ref{tab:ex05_overflow_point}. With the analytical solution as the reference, all three discretizations yield consistent predictions, with relative errors on the order of $10^{-2}$. The quadtree mesh provides the smallest error of $2.13\times10^{-2}$, whereas the polygonal and adaptive meshes exhibit similar errors of $3.39\times10^{-2}$ and $3.34\times10^{-2}$, respectively. Consistent with these results, the pressure-head contours in Fig.~\ref{fig:ex05_Pressure_head} remain smooth and physically plausible across the saturated--unsaturated transition zone, indicating that the CSFEM can robustly capture the seepage field and the associated free-surface behavior.

Tab.~\ref{tab:ex05_adaptive_efficiency} compares the computational cost of the three meshes. The quadtree mesh reduces the mesh size and runtime substantially, decreasing the number of elements from 2928 to 1904 and the CPU time from 14.8~s to 5.2~s, while preserving the accuracy level indicated by the overflow-point prediction. The adaptive mesh further improves efficiency, requiring only 767 elements and 853 DOFs and completing the analysis in 2.0~s, yet maintaining an overflow-point error comparable to that of the polygonal mesh. This demonstrates that adaptive refinement can achieve a favorable balance between accuracy and efficiency by concentrating resolution in regions of strong gradients, leading to near-equivalent solution quality at a markedly reduced computational cost.

\begin{table}[H]
\centering
\caption{Coordinates of the overflow point.}
\begin{tabular}{lcc}
\toprule
Method & Coordinate of X (m) & Relative error\\
\midrule
CSFEM (Polygonal mesh) & 0.684807 &$3.39\times10^{-2}$\\
CSFEM (Quadtree mesh) & 0.676511 &$2.13\times10^{-2}$ \\
CSFEM (Adaptive mesh) & 0.684488 &$3.34\times10^{-2}$\\
Analytical solution \cite{hornung1985evaluation} & 0.662382 & - \\
\bottomrule
\end{tabular}
\label{tab:ex05_overflow_point}
\end{table}

\begin{table}[H]
\centering
\caption{Computational cost of the CSFEM with different mesh discretizations.}
\begin{tabular}{lccc}
\toprule
Method & Elements & DOFs & CPU time (s)  \\
\midrule
CSFEM (Polygonal mesh) & 2928 & 5808 &  14.8 \\
CSFEM (Quadtree mesh) & 1904 & 2001 & 5.2  \\
CSFEM (Adaptive mesh) & 767 & 853 &  2.0  \\
\bottomrule
\end{tabular}
\label{tab:ex05_adaptive_efficiency}
\end{table}

\subsubsection{Homogeneous trapezoidal dam}
\begin{figure}[H]
    \centering
    \includegraphics[width=1.0\linewidth]{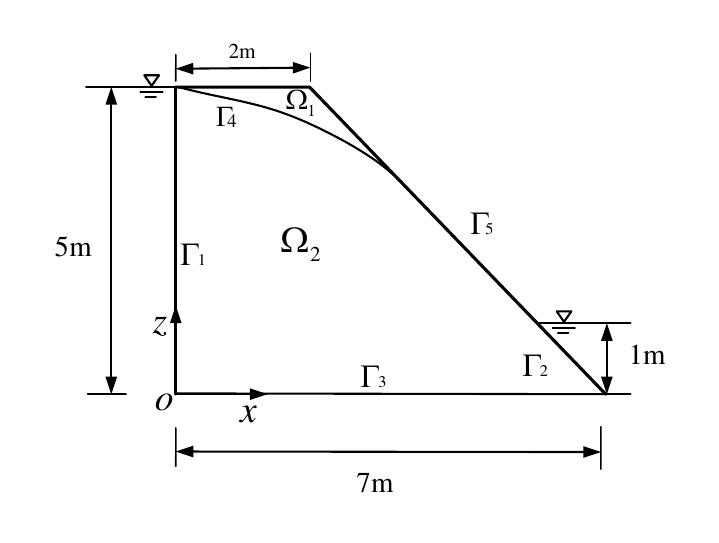}
    \caption{Geometric model and boundary conditions of a trapezoidal dam.}
    \label{fig:ex06_geo}
\end{figure}

\begin{figure}[H]
    \centering
    \includegraphics[width=0.8\linewidth]{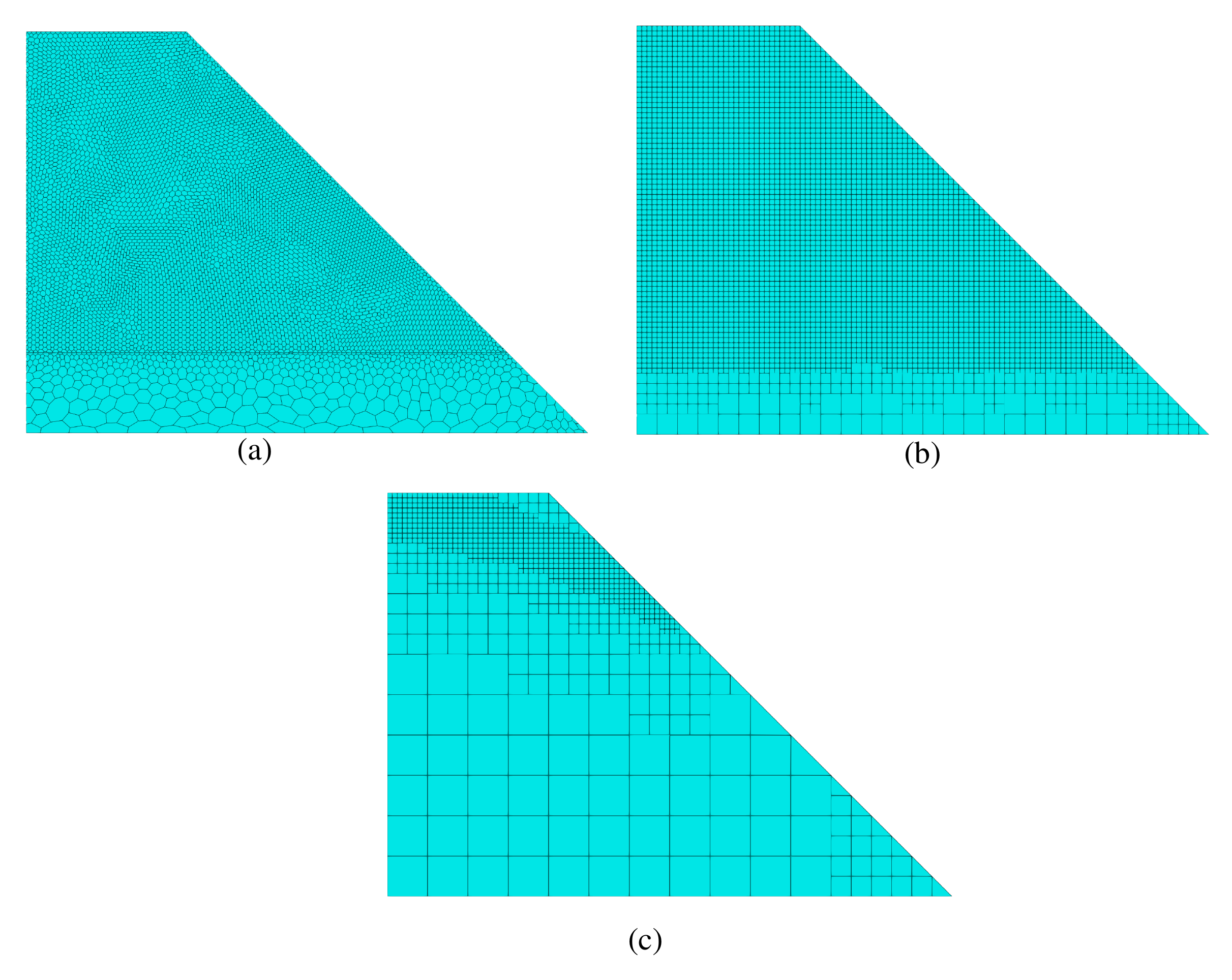}
    \caption{Mesh model of a trapezoidal dam; (a) the polygonal mesh; (b) the quadtree mesh; (c) the adaptive mesh. The adaptive mesh in (c) is the final mesh obtained by the indicator in Eqs.~\eqref{eq:band_indicator}--\eqref{eq:grad_indicator} with bulk-chasing marking (Eq.~\eqref{eq:dorfler_marking}, $\theta=0.5$) and 2:1 balanced quadtree refinement.}
    \label{fig:ex06_mesh}
\end{figure}

This example represents the cross-section of a homogeneous trapezoidal dam, with the computational model dimensions shown in Fig.~\ref{fig:ex06_geo}. The hydraulic head is 5~m on the left upstream side and 1~m on the right downstream side. The bottom is impermeable, and the permeability coefficient is 1~m/s. As shown in Fig.~\ref{fig:ex06_mesh}, the dam body is discretized using polygonal, quadtree, and adaptive meshes, where the polygonal and quadtree meshes adopt local refinement in the potential overflow region with a cell size of 0.0625~m and use a coarser cell size of 0.25~m elsewhere, while the adaptive mesh automatically concentrates refinement near the phreatic surface and potential overflow area and remains relatively coarse in the rest of the domain. 

\begin{figure}[H]
    \centering
    \includegraphics[width=1.0\linewidth]{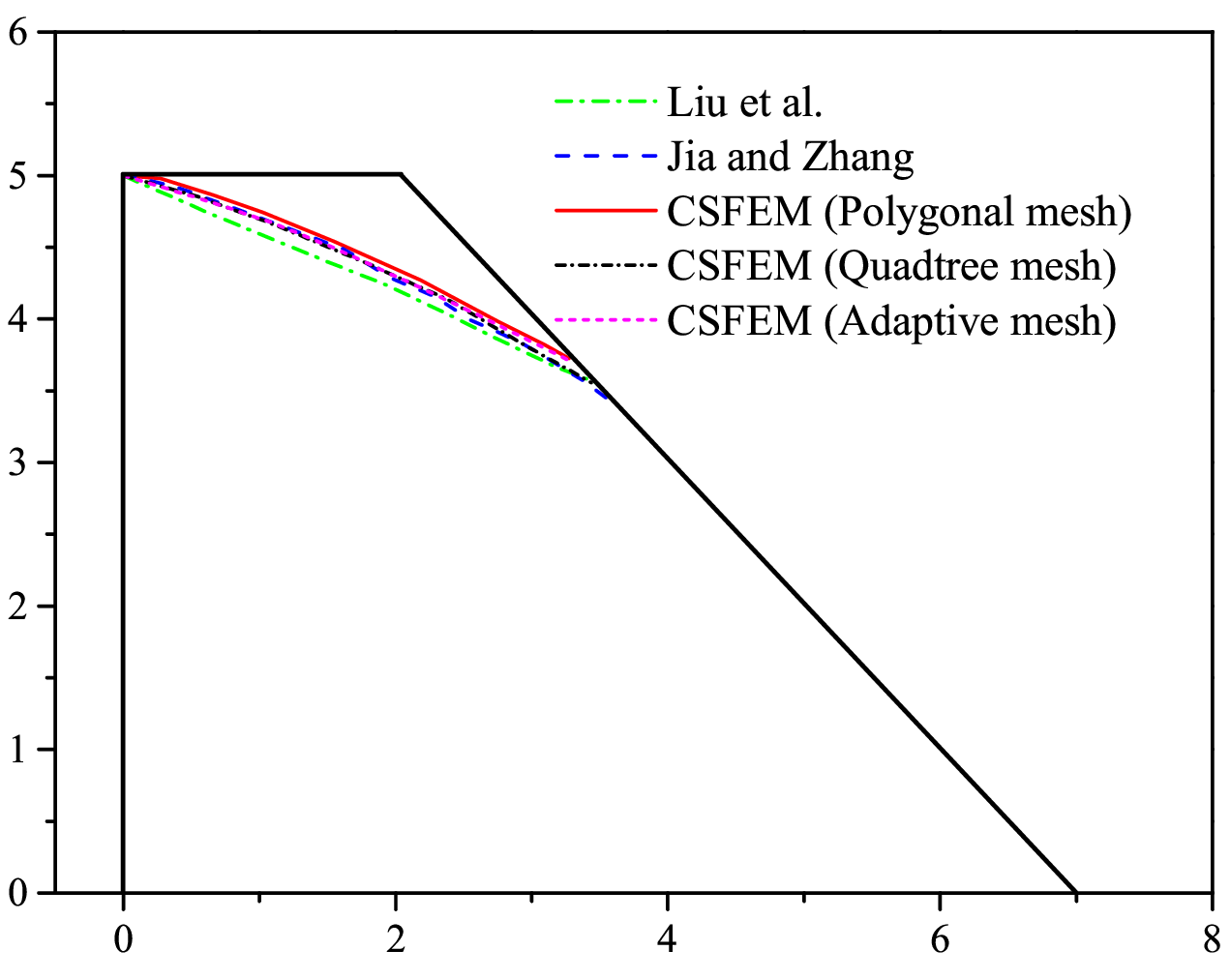}
    \caption{Comparison of free-surface positions of a trapezoidal dam.}
    \label{fig:ex06_free_surface}
\end{figure}

\begin{figure}[H]
    \centering
    \includegraphics[width=1.0\linewidth]{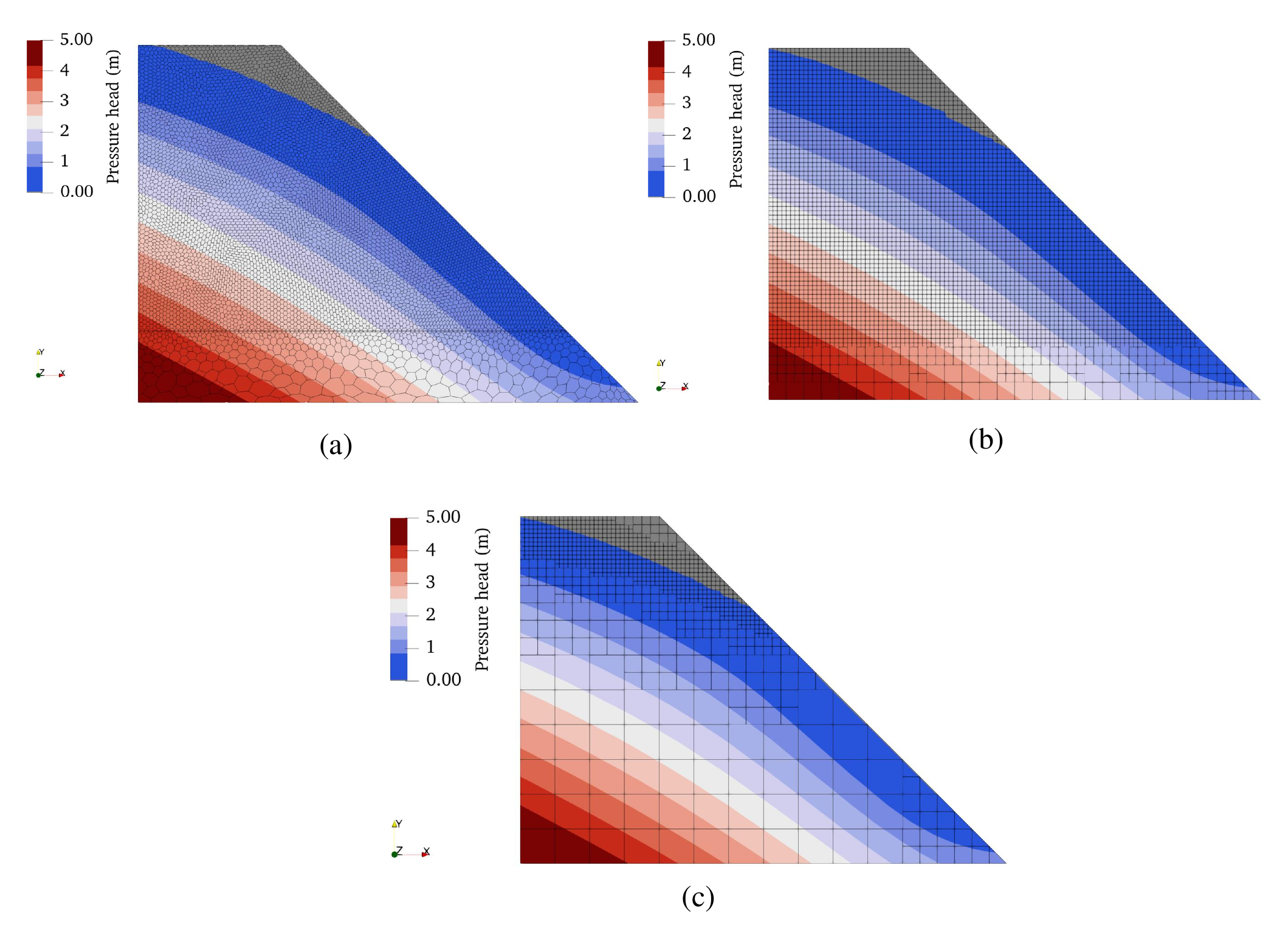}
    \caption{The pressure head distribution within a trapezoidal dam; (a) the polygonal mesh; (b) the quadtree mesh; (c) the adaptive mesh.}
    \label{fig:ex06_contour}
\end{figure}

Fig.~\ref{fig:ex06_free_surface} shows a comparison of free surface profiles predicted between the CSFEM and FEM. The CSFEM exhibits excellent consistency with the Liu et al.\citep{liu2018new} and Jia and Zheng \citep{jia2024new}. Moreover, Tab.~\ref{tab:ex06_adaptive_efficiency} reports the computational cost for different CSFEM mesh discretizations, showing that the adaptive mesh achieves the most economical computation with 767 elements and 853 DOFs and a CPU time of 2.6~s, while the polygonal and quadtree meshes require substantially more degrees of freedom and runtime. The pressure head contours presented in Fig.~\ref{fig:ex06_contour} further highlight the effectiveness of the method in capturing the seepage characteristics associated with free-surface flow.

\begin{table}[htbp]
\centering
\caption{Computational cost of the CSFEM with different mesh discretizations.}
\begin{tabular}{lccc}
\toprule
Method & Elements & DOFs & CPU time (s)  \\
\midrule
CSFEM (Polygonal mesh) & 6809 & 13520 &  72 \\
CSFEM (Quadtree mesh) & 4715 & 4828 & 30.8  \\
CSFEM (Adaptive mesh) & 767 & 853 &  2.6  \\
\bottomrule
\end{tabular}
\label{tab:ex06_adaptive_efficiency}
\end{table}

\section{Conclusions}
\label{Conclusions}

A polygonal CSFEM framework has been developed for two-dimensional seepage analyses covering steady-state, transient, and free-surface problems. By combining Wachspress interpolation with cell-based gradient smoothing, element matrices are assembled via boundary integrals only, avoiding in-element derivatives and improving robustness on distorted and locally refined meshes. Based on the numerical results presented in this work, the following conclusions are drawn.

(1) The proposed formulation passes the patch test on both polygonal and quadtree meshes, reproducing linear hydraulic-head fields with relative errors on the order of $10^{-8}$, which verifies linear completeness and numerical stability on different mesh topologies.

(2) For steady-state and transient seepage benchmarks on regular and irregular domains, including curved boundaries and complex cavities, the polygonal CSFEM solutions agree closely with FEM reference results and available analytical or literature solutions, demonstrating accuracy and robustness for practical seepage simulations. In a representative steady seepage example, the CSFEM yields smaller hydraulic-head errors than a conventional FEM at the same characteristic mesh size.

(3) For free-surface seepage in homogeneous dams, the fixed-mesh iteration coupled with the polygonal CSFEM produces smooth free-surface profiles and physically consistent pressure-head contours. Across the investigated discretizations, the computed overflow-point locations exhibit relative errors on the order of $10^{-2}$, and quadtree-based meshes generally yield slightly smoother free-surface curves due to their more regular layout and locally enhanced resolution.

(4) Quadtree refinement and solution-driven adaptivity significantly improve computational efficiency within the proposed framework. Compared with polygonal discretization, quadtree meshes reduce the degrees of freedom and runtime while preserving accuracy, and adaptive meshes further achieve comparable solution quality with markedly reduced mesh size and CPU time by concentrating refinement near the phreatic surface and other high-gradient regions.

Future work will extend the present framework to heterogeneous and anisotropic conductivity fields, develop more rigorous a posteriori error indicators for adaptive refinement, and investigate coupled hydro-mechanical applications in which evolving seepage fields interact with deformation and stability.

\section{Acknowledgements}
The Yunnan Province Xing Dian Talent Support Program (grant NO. XDYC-QNRC-2022-0764) and Yunnan Fundamental Research Projects (grant NO. 202401CF070043) provided support for this study. 

\appendix

\bibliographystyle{elsarticle-num} 
\bibliography{cas-refs}





\end{document}